\magnification=\magstep1
\input amstex
\documentstyle{amsppt} 
\NoBlackBoxes
\NoRunningHeads

\vsize=8.7truein

\def\dd{{\dot D}}
\def\R{{\Bbb R}}
\def\df{{\mathop{\ =\ }\limits^{\mathop{\text{\rm df}}}}}
\def\dist{{\mathop{\text{\rm dist}}}}
\def\var{{\mathop{\text{\rm Var}}}}
\def\cov{{\mathop{\text{\rm Cov}}}}
\def\jeden{{\mathop{\text{\bf 1}}}}
\def\prt{\partial}
\def\eps{\varepsilon}
\def\wt{\widetilde}
\def\wh{\widehat}

\def\today{\ifcase\month\or January\or February\or March\or April\or
May\or June\or July\or August\or September\or October\or November\or
December\fi\space\number\day, \number\year}

\countdef\tagno=255
\countdef\sectno=254
\tagno = 1
\sectno = 1
\def\tn#1{\the\sectno.\the\tagno \global\count#1 = \the\tagno
          \global\advance\tagno by 1 } 
\def\Tn#1{(\the\sectno.\the\tagno) \global\count#1 = \the\tagno
          \global\advance\tagno by 1  } 

\def\lbl#1{\the\sectno .\the\count#1 }
\def\Lbl#1{(\the\sectno .\the\count#1)}
\def\slbl#1#2{#1.\the\count#2 }
\def\SLbl#1#2{(#1.\the\count#2)}
\count14=3
\count12=6

\topmatter
\title
On minimal parabolic functions \\
and time-homogeneous parabolic $h$-transforms 
\endtitle
\author
Krzysztof Burdzy\\
Thomas S. Salisbury
\endauthor
\thanks
Research partially supported by NSF grant DMS-9700721, \hfill\break
and by a 
grant from NSERC 
\endthanks
\affil
University of Washington \\
York University and the Fields Institute
\endaffil
\address
Department of Mathematics, Box 354350, 
University of Washington, Seattle, WA 98195-4350
\endaddress
\email burdzy\@math.washington.edu \endemail 
\address
Department of Mathematics and Statistics, York University,
North York, Ontario, Canada M3J 1P3 
\endaddress
\email salt\@nexus.yorku.ca \endemail 
\date\today\enddate 
\abstract
Does a minimal harmonic function $h$ remain minimal 
when it is viewed as a parabolic function? 
The question is answered for a class
of long thin semi-infinite tubes $D\subset \R^d$ of variable width
and minimal harmonic functions
$h$ corresponding to the 
boundary point of $D$ ``at 
infinity.'' Suppose $f(u)$
is the width of the tube $u$ units away from its
endpoint and $f$ is a Lipschitz function. 
The answer to the question is affirmative if and only if
$\int^\infty f^3(u)du = \infty$. 
If the test fails, there exist parabolic
$h$-transforms of space-time Brownian motion
in $D$ with infinite lifetime which are not
time-homogenous.
\endabstract
\endtopmatter

\document
\subheading{1. Introduction and main results}
We want to compare the parabolic Martin boundary of
a domain in $\R^d$
with its Martin boundary, both topologically and
probabilistically. In many cases, the two boundaries are
related in a very simple way. This provides a complete
description
of the parabolic Martin boundary in those cases (quite
many) when the Martin boundary is known.
We plan to present a detailed discussion
of this general problem in a separate publication.
This paper is devoted to a 
narrower 
aspect of the relationship between the two boundaries.
We will start with a very informal discussion
of a special case which motivated our study.
The concepts
of the usual and parabolic Martin boundary will be reviewed 
in a rigorous way later in the introduction. 
The basic ideas of the classical potential theory and
Brownian motion may be found in Doob (1984).

Consider a strip $D = \{(x^1, x^2)\in \R^2: |x^2| < 1\} $.
Let $X_t$ be a Brownian motion starting from $(0,0)$.
Then $\dot X_t =(X_t, -t)$ is a space-time Brownian motion
starting from $(0,0,0)$. First fix some $s>0$, a point 
$z \in \prt D$ 
and a sequence of points 
$\{z_k\}$ in $D$ converging to $z$ 
as $k\to \infty$. Condition $\dot X$ to be at 
$(z_k, -s)$ at time $s$ and to 
not leave $D\times \R$
before time $s$. Then let $k$ go to infinity. The conditioned
processes converge in distribution to a process whose
first 
coordinate 
is a Brownian motion conditioned to
exit $D$ through 
$z$ at time $s$. The lifetime 
of this process is finite. This conditioned
space-time Brownian motion is not time-homogeneous,
i.e., its transition probabilities
$P(\dot X_u \in (dy,-du) \mid \dot  X_t \in (dx, -dt) )$
depend not only on $u-t$, but on the values of $t$ and $u$ 
as well.

Next suppose that $c>0$ is a constant and
consider $\dot X$ conditioned to be at $(ck, 0, -k)$
at time $k$ and 
to 
not leave $D\times \R$
before time $k$. In the limit, 
as 
$k\to\infty$,
we obtain a process whose 
spatial component escapes ``to $+\infty$'' within $D$
at rate $c$. The first coordinate 
of the space process is a one-dimensional Brownian motion with
drift $c$. This conditioned space-time Brownian motion 
is time-homogeneous and its lifetime is infinite.

The domain in our example, a strip, seems to be typical
and we would expect that many domains have the property
stated in the following problem.

\proclaim{\Tn{11} Problem} 
Find necessary and sufficient 
conditions, of a geometric nature in $D$, such that for every 
minimal
parabolic function $h$ in $\dd$, the corresponding
$h$-transform of the space-time Brownian motion is time
homogeneous if and only if its lifetime is a.s. infinite.
\endproclaim

Another source of 
motivation may be explained in
purely analytic language. Recall the domain 
of 
our
first example,
$D = \{(x^1, x^2)\in \R^2: |x^2| < 1\} $.
Consider a minimal positive harmonic
function $h(x)$, $x\in D$. Let $g(x,t) = h(x)$ for
all $x \in D$ and $t \in \R$. Evidently, $g$ is a parabolic
function, and we may therefore identify every harmonic function with 
a parabolic function.
Since $h$ is minimal harmonic,
it corresponds to a minimal Martin boundary point $y$ of $D$.
Suppose that $y$ is also a Euclidean boundary point, say,
$y= (1,1)$. Then 
$g$  
is not minimal as a parabolic function,
i.e., it is a mixture of different parabolic functions.
An easy probabilistic justification 
can be based on the 
fact that Brownian motion conditioned
by $h$ has a random 
lifetime. Thus 
the space-time Brownian 
motion conditioned by 
$g$ 
is a mixture of processes
conditioned to exit $D$ through $y$ at different times $s$,
i.e., a mixture of 
$g_s$-transforms for different parabolic functions $g_s$. 
However, if $y$ is the point at ``$+\infty$'' then 
$g$ 
is minimal in the space of parabolic functions. 
While not completely obvious, this is simple to show directly, and
also follows from our main result, Theorem \lbl{14} below.  
Our informal discussion
suggests that in many domains, a minimal harmonic
function is also minimal in the space of parabolic
functions if and only if it corresponds to a ``point
at infinity.'' We propose the following problem.

\proclaim{\Tn{13} Problem}Determine which minimal harmonic
functions are minimal in the space of parabolic functions.
\endproclaim

We are not able to give a complete answer to either of the 
two problems but we hope that our main result, Theorem \lbl{14}
below, will shed light on both. 

We proceed with a rigorous presentation of our results.
We start with a review of basic definitions
and facts concerning Martin boundaries and conditioned
Brownian motion.
Let $D$ be a Euclidean domain, that is, an open connected
subset of  $\R^d$ for some $d \geq 2$. We will consider the
domain $\dd \df D \times (-\infty, 0)\subset \R^{d+1}$. Let
$G(x,y)=G_D(x,y)$ and $\dot G(u,v)=\dot G_{\dd}(u,v)$ be the Green 
functions for 
$(1/2)\Delta$ on $D$ and for the heat operator $(1/2)\Delta
-\partial/\partial t$ on $\dd$ where $\Delta$ is the Laplace
operator (see Doob (1984) 1.VII.1 and
1.XVII.4). Thus $G : D \times D \to (0, \infty]$ and $\dot G :
\dd \times \dd \to [0,\infty] $.  For $u=(x,s)\in \dd$ and
$v=(y,s-t)\in\dd$ we have that
$$
\dot G(u,v) = \cases p_t(x,y),&\text{for $t>0$}\\
0,&\text{for $s<t\leq 0$},\endcases
$$
where $p_t=p_t^D$ is the heat 
kernel on $D$ (that is, the transition function for Brownian
motion killed upon leaving $D$). 
Note that this formula can also be used to define $\dot G((x,s),v)$ 
when $s=0$. 
A function $h:D\to[0,\infty)$ is {\sl harmonic} if $\Delta h=0$ on
$D$. A function $g:\dd\to[0,\infty)$ is {\sl parabolic} if it solves
the heat equation
$$
\frac{\partial g}{\partial t}=\frac{1}{2}\Delta_x g
$$ 
in $\dd$. 
In this case, it is {\sl superparabolic} as well. That is,
$$
g(x,s)\ge\int g(y,s-t)p_t(x,y) dy
$$ 
for every $(x,s)\in\dd$ and $t>0$. 
We may extend $g$ by letting 
$$
g(x,0)\df \lim_{t\downarrow 0}\int g(y,-t)p_t(x,y)dy
$$
(the limit is easily seen to be monotone).
We say that $g$ is {\sl admissible} if $g(x_0,0)<\infty$. 

Now recall the definitions of the Martin boundary in the
elliptic and parabolic contexts (Doob (1984) 1.XII.3 and
1.XIX.3). Fix some $x_0 \in D$ and let
$$
K(x,y) \df \frac{G(x,y)}{G(x_0,y)}
$$
for $x,y \in D$. Then,
up to homeomorphism there is a unique  metrizable
compactification $D^M$ of $D$ such that
\roster
\item"(i)" the function $K(\,  \cdot \, , \, \cdot \, )$ may be
extended continuously to $D\times (D^M \setminus \{x_0\} ) $;
\item"(ii)" $K(\,\cdot\,,  x) \equiv K(\,\cdot\,,  y)$ if
and only if $x = y$.
\endroster
The set $\prt^M D \df D^M \setminus D$ is called the 
{\sl Martin boundary} of $D$. 
For $z\in\prt^M D$ and $y_k\in D$, we have $y_k\to z$ if and only if
$K(x,y_k)\to K(x,z)$ for every $x\in D$. A harmonic function 
$h>0$ is said to be minimal if, whenever $h'>0$ is harmonic, and
$h'\le h$, it follows that $h'=ch$ for some constant $c$. A point
$z\in\prt^M D$ is said to be minimal if $K(\,\cdot\, ,z)$ is minimal. 
For every $h>0$ harmonic, there is a unique measure $\mu$,
concentrated on the set $\prt^M_0 D$
of minimal points of $\prt^M D$, such that
$$
h(x)=\int_{\prt^M_0 D} K(x,z)\mu(dz),
$$
for every $x\in D$ (See Doob (1984) 1.XII.9). 

Now 
define $\dot K$ on $\dd\times\dd$ by
$$
\align
\dot K((x,s), (y,t)) &\df \frac{\dot G ((x,s),(y,t))} {\dot
G ((x_0, 0), (y,t))}\\
&=\cases
p_{s-t}(x,y)/p_{-t}(x_0,y),&t<s<0\\
0,&s\le t<0.
\endcases 
\endalign 
$$
Then up to homeomorphism, there
is a unique  metrizable compactification $\dd^M$ of $\dd$ with
the following properties:
\roster
\item"(i)" the function $\dot K$ has
an extension to $\dd \times \dd^M$ such that for each
$(x,s)\in\dd$, the function $\dot K ((x,s),\, \cdot\, )$ is
finite valued and continuous on $\dd^M \setminus \{(x,s)\} $;
\item"(ii)" $\dot K( \, \cdot\,,u) = \dot K(  \,
\cdot\,,v)$ if and only if $u=v$.
\endroster
\noindent We call $u$ the {\sl pole} of $\dot K( \, \cdot\,,u)$. 
We write $\prt^M\dd \df \dd^M \setminus\dd$ and call it the
{\sl Martin boundary} of 
$\dd$ (or the parabolic Martin boundary of
$D$).  
We have again that, for $z\in\prt^M\dd$ and $(y_k,t_k)\in \dd$,
$(y_k,t_k)\to z$ if and only if $\dot K((x,t),(y_k,t_k))\to\dot
K((x,t),z)$ for every $(x,t)\in\dd$. Every $\dot K(\,\cdot\,,z)$ is
admissible (see 1.XIX.3.1 of Doob (1984). 

We denote by $\dot 0$ the unique point of
$\prt^M\dd$ for which $K(\, \cdot\,,\dot 0)\equiv 0$.
It is unique by (ii) and exists as the limit of
some subsequence of $(x_0,1/n)$.  
A point $z\in\prt^M\dd$ is {\sl minimal} if $\dot K(\,\cdot\,,z)$ is  
minimal as a parabolic function, and $\dot K((x_0,0),z)=1$.
The set of minimal points is denoted $\prt^M_0\dd$. The integral 
representation of admissible parabolic functions as 
$$
g(x,t)=\int_{\prt^M_0\dd}\dot K((x,t),z)\mu(dz)
$$
is entirely analogous to that of 
the harmonic setting (See Doob (1984)
1.XIX.7).

Let $(\Omega, \Cal F)$ be a measurable space with 
$X:\Omega \times [0,\infty) \to \R^d\cup\{\delta\}$ 
a stochastic process. 
We use the notation $X_t$ and $X(t)$ interchangeably.
$P^x$ is a probability measure under which $X$ is a standard
$d$-dimensional Brownian motion started from $x$, and killed
upon leaving $D$. We write $E^x$ for the corresponding expectation. 
In particular,
$\delta$ is a cemetery point adjoined to $D$, $X$ is
continuous on a random time interval $[0,\zeta)$, and
$X_t=\delta$ for $t\ge\zeta$.

Let $\tau_t=\tau_0-t$ be a process measuring absolute time, and
write $\dot X_t=(X_t,\tau_t)$.
By enlarging $\Omega$ if necessary, we may suppose 
that for each $s\le 0$, there are probability measures $ P^{x,s}$
under which
$X$ has the same law as under $P^x$, and $\tau_0=s$. That is,
$\{\dot X_t, t\geq 0\}$ is a space-time Brownian motion
starting from $(x,s)$. 

If $h:D\to (0,\infty]$ is a superharmonic function then
$$
p^h_t(x,y)\df
\frac{h(y)p_t(x,y)}{h(x)}
$$
is the transition function 
of a Markov process $X^h$, called 
an $h$-transform, or conditioned
Brownian motion. 
We write $P^x_h$ and $E^x_h$ for the corresponding probability
measure, and its expectations. 
By convention, $h$ is taken to vanish at
$\delta$. If $x\in D^M$, $x\ne x_0$ then we write $X^x$ for
$X^{K(\, \cdot \,,x)}$. 
If $h=\int_{\prt^M_0}K(\,\cdot\, ,z)\mu(dz)$ is harmonic, then
$$
P^x_h=\frac{1}{h(x)}\int_{\prt^M_0} K(x,z) P^x_z\,\mu(dz).
$$ 
The paths of $X^h$ converge a.s. to points of the minimal Martin
boundary, at their lifetimes (see Doob (1984) 3.III.1, or
section 7.2 of Pinsky (1995)). 

Similarly, if $g:\dd\to [0,\infty]$ is a  superparabolic
function, then
$$
\dot p^g_t((x,s),(y,s-t))\df
\frac{g(y,s-t)p_t(x,y)}{g(x,s)}
$$
is the transition
function for a Markov process $\dot X^g$ taking values in
$\dd\cup\{\delta\}$ (actually in
$\{\delta\}\cup\{u\in\dd;\,g(u)>0\}$) that we call a
conditioned space-time Brownian motion. We will use
$P_g^{x,s}$ to denote a probability measure under which $\dot
X^g$ has this transition function and starts from $(x,s)$. We
write $X^g$ for the spatial component of $\dot X^g$ (with
$X^g_t=\delta$ for $t\ge\zeta$), and note that
$$
\dot X^g_t =\cases
(X^g_t,\tau_t)\in\dd,&\text{for }t<\zeta\\ 
\delta &\text{for } t\ge\zeta.\endcases
$$
We will also refer
to $X^g$ as an $g$-transform. This abuse should cause no
confusion, as it is easy to check that if $h$ is superharmonic
and we define a superparabolic function $g$ by $g(x,t)=h(x)$
then $X^h=X^g$. If $u\in\dd^M$ then we write $\dot X^u$,
$X^u$, $P^{x,s}_u$ instead of $\dot X^{\dot K(\, \cdot
\,,u)}$, etc. Strictly speaking, the above formulae hold
under $P_g^{x,s}$ only for $s<0$, but by taking $X^g_0=x$
under $P_g^{x,0}$, we obtain extensions valid for $s=0$ as
well, provided $g$ is admissible. 
If $g$ is actually parabolic, then each $g$-process
approaches the one-point boundary of $\dd$ at its lifetime
$\zeta$ 
(Doob (1984) 2.X.12), in other words, it eventually leaves every
compact subset of $\dot D$. In the Martin topology, the paths of
$\dot X$ converge at their lifetimes, to points of the minimal
parabolic Martin boundary, and the measures $P^{x,s}_g$ can be
represented in terms of the $P^{x,s}_u$, for $u\in\prt^M_0\dd$, just
as in the harmonic setting. 

For $(x^1, x^2, \dots , x^d) \in \R^d$ let
$\wt x = (x^1,x^2, \dots, x^{d-1} )$.
We will restrict our attention to ``tubes'' with variable width.
For a non-negative function $f: \R \to \R$, let
$$D_f \df \{ x\in R^d : |\wt x| < f(x^d) \} .$$
We will always assume that $f$ is
strictly positive on $(a,b)$ for some 
$-\infty \leq a < b \leq \infty$ and equal to 0 on 
$(-\infty,a] \cup[b,\infty)$.
We will focus on domains $D_f$ corresponding to functions $f$
which are Lipschitz on $(a,b)$ (the function may have a jump
at $a$ or $b$).
If $f$ is Lipschitz and $b=\infty$, then each sequence $x_k$ of points in
$D_f$ such that $x_k^d \to \infty$ converges in the Martin
topology to a point (the same for all such sequences) which we
will denote as $\infty$. 
The proof of this claim is easy --- it
may be based on the boundary Harnack principle. 
The same result should be true 
for all functions $f$ (not necessarily
Lipschitz) but we do not see an obvious argument. An analogous
remark applies to $-\infty$. 
Any positive harmonic function $h$ corresponding to
$\infty \in \prt^MD_f$
vanishes on $\{x\in \prt D_f: x^d <b\}$ and, moreover,
$h(x) \to 0$ when $x^d \to -\infty$.

Let $\Lambda_s = \{x\in D_f : x^d =s \} $.
The stopping time $\inf \{t>0: X_t \in A\}$ will be denoted
$T(A)$.
We write $\tau(A)$ for the absolute time $\tau_{T(A)}=\tau_0-T(A)$. 

Recall that a harmonic function $h$ is identified with a parabolic
function by letting $h(x,t) =h(x)$.

\proclaim{\Tn{14} Theorem} Suppose that $b=\infty$
and $f$ is a function which is Lipschitz on $(a,b)$ and such that 
$$\limsup_{v\to\infty} f(v) <\infty$$
and 
$$\int_u^\infty f(v) dv =\infty
\tag\tn{10}$$ 
for
all $u<\infty$. Let $h$ be the
minimal harmonic function corresponding to $\infty \in 
\prt^M_0 D_f$. 
Fix some $x_0 \in D_f$. 

\medskip\noindent (i) Suppose that either
\roster
\item"(a)" $\int_a^\infty f^3(v)dv < \infty$ or
\item"(b)" the Lipschitz
constant of $f$ is sufficiently small (it will suffice
to assume that it is less than the $\lambda$ in 
(iv) of Theorem \lbl{12}) and 
$\int_u^\infty f^3(v)dv < \infty$ for some $u<\infty$.
\endroster
Each one of assumptions (a) or (b) implies (A)-(D) below.
\roster
\item"(A)" For some function 
$g: (a,\infty) \to (-\infty, 0]$ 
with $\lim_{u\to\infty}g(u)=-\infty$,
we have the following.
For each $s\in \R$ there is a minimal point $z_s \in \prt^M_0 \dd_f$,
which is the limit of all sequences 
$(x_k, (g(x_k^d)-s_k)\land 0)$ with  $x_k^d \to \infty$ and
$s_k\to s$. 
\item"(B)" 
If $s_1 \ne s_2$ then $z_{s_1} \ne z_{s_2}$. 
\item"(C)" Let $h_s$ denote a minimal parabolic function
with pole at 
$z_s$. Then $h = \int_\R h_s \mu(ds)$ for some
measure $\mu$ which charges all non-degenerate intervals.
In particular, $h$ is not
minimal in the space of parabolic functions on $\dot D_f$.
\item"(D)" Let $s\in\R$ and $(x,t)\in\dd$. The process
$\dot X$ is not time-homogeneous under $P_{z_s}^{x,t}$. In fact,
$g(u)-\tau(\Lambda_u)\to s$ as $u\to\infty$ $P_{z_s}^{x,t}$-a.s.
Hence,
$\lim_{u\to\infty} (T(\Lambda_u)+g(u))$ exists $P_h^x$-a.s. 
\endroster
(ii) If $\int_u^\infty f^3(v)dv = \infty$ for all $u<\infty$
then $h$ is 
minimal in the space of parabolic functions on $\dot D_f$.
\endproclaim

\remark{\Tn{15} Remarks} 

The lifetime of Brownian motion conditioned by $h$ is
infinite if and only if $\int_u^\infty f(v) dv =\infty$ for
all $u<\infty$, according to 
Theorem \lbl{12} 
below. If this condition is not 
satisfied, the function $h$ is not minimal as a parabolic function 
(see the discussion preceding 
Problem \lbl{13}). 
\endremark

The proof of 
Theorem \lbl{14}  hinges on 
estimates of the variance
of $h$-path lifetimes. Since the estimates may have some 
independent interest, we state them as 
Theorem \lbl{12} below.

Several authors 
have 
addressed the problem of 
when, 
given a domain $D\subset \R^d$,
there is a constant $c=c(D)<\infty$
such that for any $x\in D$ and any positive harmonic
function $h$ in $D$ we have $E^x_h \zeta < c$. 
The pioneering work 
was 
done by
Cranston and McConnell (1983) and Cranston (1985).
The existence of the finite upper bound $c$ is known for a wide class of
domains; 
see, e.g., Ba\~nuelos and Davis (1992) or
Bass and Burdzy (1992) and references therein. Higher
moments of $h$-path lifetimes have been studied by
Davis (1988), Davis and Zhang (1994) and Zhang (1996).

Chris Rogers has pointed out to us that a 
related equivalence, between non-minimality and the variance of
hitting times, has been established in the context of
one-dimensional diffusions. There, the speed measure and coupling
can be used to give a simple proof. See Rogers (1988), which 
synthesizes earlier work of Fristedt and Orey (1978), K\"uchler and
Lunze (1980), and R\"osler (1979). 

Recall that we are concerned with functions $f$ which are 
strictly positive and Lipschitz on $(a,b)$ and equal to 0 on 
$(-\infty,a] \cup[b,\infty)$.
Our next result holds for all functions $f$ which
are Lipschitz on $(a,b)$. However,
in order to simplify the notation we will prove it only in the case
when $f$ is Lipschitz with the constant equal to $1$, i.e.,
from now on we will assume that
$|f(u) - f(v)| \leq |u-v|$ for $u,v \in (a,b)$.
Fix some $s_0 \in (a,b)$ and define $s_k$ inductively by
$s_{k+1} = s_k + f(s_k)/2$ for $k \geq 0$ and
$s_{k-1} = s_k - f(s_k)/2$ for $k \leq 0$.
If $s_k\geq b$ for some $k$ then we redefine $s_j$ for
$j\geq k$ and we let $s_j=b$ for all
$j\geq k$. A similar remark applies to the case when $s_k\leq a$.
Note that it may happen that $s_k < b$ for all $k>0$
and/or $s_k > a $ for all $k<0$. However, we always have
$\lim _{k\to \infty} s_k = b $ and
$\lim _{k\to -\infty} s_k = a $. Let
$k_f = \inf\{k: s_k = b\}$ and recall that
$\Lambda_{s_k} = \{ x\in D_f: x^d = s_k\}$.
Let $D_j$ be the
component of $D_f \setminus \Lambda_{s_j}$ which contains points $x$
with $x^d < s_j$.

\proclaim{\Tn{12} Theorem} Let $h$ be a positive harmonic
function in $D_f$ which vanishes on $\{x\in \prt D: x^d <b\}$.
If $b = \infty$ then $h$ corresponds to $\infty \in 
\prt^M_0 D_f$. 
In the following statements, $x$ ranges over the elements of $D_f$ 
with $x^d < b - f(b-)$ (here
$\infty-\infty=\infty$).

\roster
\item"(i)" For some $c_1,c_2 \in (0,\infty)$,
$$c_1 \int_{x^d}^b f(v) dv \leq E^x_h \zeta \leq
c_2 \int_{x^d}^b f(v) dv. \tag\tn{17}$$ 

\item"(ii)" If $\int_{x^d}^b f(v) dv = \infty $
then $\zeta = \infty$ $P^x_h$-a.s.

\item"(iii)" If $\zeta <\infty $ $P^x_h$-a.s. then
for some $c_3,c_4 \in (0,\infty)$,
$$c_3 \int_{x^d}^b f^3(v) dv \leq \var^x_h  \zeta \leq 
c_4 \int_a^b f^3(v) dv .\tag\tn{18}$$ 

\item"(iv)" 
There exists $\lambda>0$ such that if the Lipschitz constant of $f$ is
less than $\lambda$ then
$$\var^x_h \zeta \leq c_5 \int_{x^d}^b f^3(v) dv .\tag\tn{37}$$

\item"(v)" If $\int_{x^d}^b f^3(v) dv = \infty$ then for each
$c_6 <\infty$ and $c_7 >0$ there is 
a  
$k_0<\infty$  such that for all $k>k_0$ and $u\in \R$,
$$P^x_h (T(\Lambda_{s_k}) \in (u, u+c_6)) < c_7.$$
\endroster
\endproclaim

\remark{\Tn{38} Remarks} 

(i) The constants $c_j$ in Theorem \lbl{12}
depend only on the dimension $d$ and the
Lipschitz constant of $f$. However, the proof will be given only in the
case when the Lipschitz constant of $f$ is equal to 1 so all the
constants in Section 2 will depend only on the dimension $d$.

(ii) The bound \Lbl{37} holds for $d\geq4$ without any 
assumptions on the value of the Lipschitz constant of $f$
but it does not hold without such an assumption
for $d < 4$. We are not going to prove 
the latter. 
It essentially
follows from 
a 
theorem of Davis and Zhang (1994).

(iii) We can give a meaning to \Lbl{18} and \Lbl{37} even if
$\zeta =\infty $ $P^x_h$-a.s. Note that in such a case we
necessarily have $b=\infty$ (see \Lbl{17}).
For all $k<\infty$ and $x\in D_f$ such
that $x^d < s_k$,
$$\var^x_h T(\Lambda_{s_k})  < c_4 \int_a^b f^3(v) dv $$
with the same constant $c_4$ as in \Lbl{18}.
This and the analogous modification of \Lbl{37} 
can be proved by applying the theorem
to the function 
$\wt f(v) \df f(v) \jeden_{(-\infty, s_k)}(v)$.

(iv) In the two-dimensional case, part (i) of Theorem \lbl{12} is due to
Xu (1990). This was generalized in Ba\~nuelos and Davis (1992).

(v) Suppose that $d=2$, the Lipschitz constant of $f$ is small
and let $\rho$ be the supremum of areas
of discs contained in $D_f$. Then \Lbl{17} and \Lbl{37} imply
that $\var^x_h \zeta \leq c_1 \rho E^x_h \zeta$. Davis (1988)
discovered this inequality and proved that it holds for all simply
connected planar domains $D$ provided $h$ is a
minimal positive 
harmonic function or a Green function.

\endremark

We would like to thank Rodrigo Ba\~nuelos, Rich Bass
and Burgess Davis for 
some very 
useful discussions
of $h$-path lifetimes.


\sectno = 2
\tagno = 1

\subheading{2. Moments of $h$-transform 
lifetimes} 
This section contains the proof of 
Theorem \slbl{1}{12}. 
We start with a short review of some useful facts about
$h$-processes.   The proofs 
may be found in Doob (1984) and Meyer, Smythe and Walsh (1972).

Let $D \subset \R^d$ be a Greenian domain and $h$ be a positive 
superharmonic function in $D$.   
Suppose that $M$ is a closed subset of $D$ and let
$L = \sup \{ t < \zeta : X _t \in M \}$
be the last exit time from $M$.  Let

\itemitem {} $Y_1 (t) = X (t), \quad t \in (0, T (M))$,

\itemitem {} $Y_2 (t) = X (T (M) + t), \quad t \in (0, \zeta-T (M))$,

\itemitem {} $Y_3 (t) = X (t), \quad t \in (0, L)$,

\itemitem {} $Y_4 (t) = X (L + t), \quad t \in (0, \zeta-L)$,

\itemitem {} $Y_5 (t) = X (\zeta-t), \quad t \in (0, \zeta)$.

\noindent Under $P_h^x$, each process $Y_k$ is an 
$h_k$-transform 
in a domain 
$D_k$, where
$D_1 = D_4 = D \backslash M$ and $ D_2 = D_3 = D_5 = D$.
Moreover, $h_1 = h_2 = h$. The function $h_3$ is a potential 
supported by $\partial M$. 
The function $h_4$ is harmonic and
has the boundary values $0$ on $\partial M$ and the same 
boundary values as $h$ on $\partial D \backslash \partial M$.
The function $h_5$ is the Green function 
$G_D (x, \,\cdot\,)$ 
if $x \in D$ or a 
harmonic function with a pole at $x$ if $x \in \partial D$.

If $\mu (dy)$ is the $P^x$-distribution of 
$X(T(M))$ 
then the $P^x_h$-distribution of this random variable is $\mu 
(dy) h (y) / h (x)$.

\proclaim{\Tn{80} Lemma} (Brownian scaling) Suppose $h$ is a 
positive superharmonic
function in a domain $D\subset \R^d$ and $x\in D^M$. For a fixed
$a\in(0,\infty)$ let
$$\align
D_a &\df \{y\in\R^d: y/a \in D\}, \\
h_a(y) &\df h(y/a) \quad \text{for} \ \ y\in D_a , \\
x_a &\df ax, \\
X^a_t &\df a X_{t/a^2} \quad \text{for} \ \ t\geq 0 .
\endalign $$
If $X$ has the 
distribution $P^x_h$, then $X^a$ has the
distribution $P^{x_a}_{h_a}$. 
\endproclaim
\demo{Proof} The lemma follows immediately from the scaling properties of
Brownian motion and superharmonic functions. \qed
\enddemo

A domain 
$D\subset\R^d$, 
$d\geq2$,
is called a {\sl Lipschitz domain} if for every $x\in\prt D$
there is a neighborhood $U_x$ of $x$, an orthonormal
coordinate system $CS_x$ and a
Lipschitz function $f_x:\R^{d-1}\to \R$ 
with constant $\lambda$ (independent of $x$) such
that $\prt D\cap U_x$ is a part of the graph of $f_x$ in $CS_x$.
Note also that the index on any constant $c_1$, $c_2$,\dots is local
in nature. That is, new results or sections of proofs will start
numbering their constants with $c_1$ as well. 

\proclaim{\Tn{81} Lemma} (Boundary Harnack principle) 

\noindent 
(a) 
Suppose
$f:\R^{d-1}\to\R$ is a Lipschitz function with constant $\lambda >0$,
$|f(x)| \leq 1$ for all $x\in \R^{d-1}
$,  
and let
$$\align
D &= \{x\in\R^d : |\wt x| < 1, f(\wt x ) < x^d < 2\} , \\
D_1 &= \{ x \in D: |\wt x| < 1/2, x^d < 3/2 \}.
\endalign $$
There exists $c_1>0$ which depends on $\lambda$ but otherwise does not
depend on $f$ such that for all $x,y\in D_1$ and all positive harmonic
functions $g,h$ in $D$ which vanish continuously on
$\{z\in \prt D: z^d = f(\wt z)\}$ we have
$$ \frac{g(x)}{g(y)} \geq c_1 \frac {h(x)}{h(y)}. $$

\noindent 
(b) 
Suppose $D$ is a Lipschitz domain, $Q$ is a compact
set and $A$ is an open set such that $Q\cap\overline D \subset A$.
There exists $c_2>0$ such that for all $x,y\in Q\cap D$ and
all positive harmonic functions $g,h$ in $D$ which vanish 
continuously on $\prt D\cap A$ we have
$$ \frac{g(x)}{g(y)} \geq c_2 \frac {h(x)}{h(y)}. \qed$$
\endproclaim

For the first proofs of the boundary Harnack principle, see
Ancona (1978), Dahlberg (1977) and Wu (1978). Stronger versions of the
result may be found in Bass and Burdzy (1991) or Ba\~nuelos,
Bass and Burdzy (1991).

Part (a) of Lemma \lbl{81} 
holds (with the same $c_1$) in domains which may be obtained from
$D$ by scaling. 

When applying the boundary Harnack principle we will sometimes leave it
to the reader to find the right choice of $D$ and $D_1$
or $D$, $A$ and $Q$.

\proclaim{\Tn{26} Lemma}
Suppose $D$ is a domain, $D_1$ is a Lipschitz subdomain
of $D$, $Q$ is a compact set, $A$ is an open
set such that $Q\cap \overline D\subset A$,
$A\cap D\subset D_1$, and $M$ is a Borel subset of 
$D\setminus A$. Assume that $h$ is a positive superharmonic 
function  in $D$ which
vanishes on $\prt D\cap A$ and is harmonic in $D_1$. Then
$$
P^x_h (T(M)<\infty )
\leq c_1
P^y_h (T(M)<\infty )
$$
for all $x,y\in Q\cap D$. The constant $c_1$ depends only 
on $D_1,Q$ and $A$.
\endproclaim
\demo{Proof}
The function
$$x \to
E^x [T(M)< T(\prt D), h(X(T(M))) ]
$$
is positive and harmonic in $A\cap D$ and the same is true for 
$x\to h(x)$. Let $D_2$ be a Lipschitz subdomain of $A\cap D$ which
contains $Q$. By the boundary Harnack principle 
\Lbl{81}(b),
applied
in $D_2$,
$$\align
P^x_h (T(M)<\infty )
& = \frac{1}{h(x)}
E^x [T(M)< T(\prt D), h(X(T(M))) ] \\
&\leq c_2
\frac{1}{h(y)}
E^y [T(M)< T(\prt D), h(X(T(M))) ] \\
&= c_2
P^y_h (T(M)<\infty ) . \qed
\endalign$$ 
\enddemo

\proclaim{\Tn{27} Lemma} Suppose $D$ is a domain and for each $k=1,2$,
\roster
\item"(i)" $D_k$ is a subdomain of $ D$,
\item"(ii)" $A_k \df \prt D_k \cap D$,
\item"(iii)" $V_k$ is an open set and $Q_k$ is a compact set
such that $Q_k \cap \overline D\subset V_k$
and $\overline V_k \cap D \subset D_k$,
\item"(iv)" $(D_1 \cup V_1) \cap (D_2 \cup V_2) = \emptyset$,
\item"(v)" there is 
a 
$c_k>0$ such that for all
$x,y\in Q_k\cap D$ and all positive harmonic functions $f,g$
in $D_k$ which vanish on $ V_k \cap \prt D$ we have
$$
\frac{f(x)}{f(y)} \geq c_k \frac{g(x)}{g(y)} .
\tag\tn{28}$$
\endroster
Assume that $x_1, x_2 \in \overline{Q_1\cap D}$ and $h_1,h_2$ are 
positive superharmonic functions in $D$ which vanish
continuously on $\prt D\setminus V_2$ and are harmonic in
$D\setminus Q_2$.
Let $T_1 \df T(A_1)$ and let $T_2$ be the last exit
time from $A_2$. The distributions of 
$\{X_t, t\in[T_1, T_2]\}$ under
$P^{x_1}_{h_1}$ and $P^{x_2}_{h_2}$
are mutually absolutely continuous and their Radon-Nikodym derivative is
bounded below by $c_1 c_2$.
\endproclaim
\demo{Proof} We will consider only the case when
$x_k \in Q_1 \cap D$ and $h_k (\,\cdot\,) = G_D(\,\cdot\, ,y_k)$
for some $y_k \in Q_2 \cap D$. Other points $x_k$
and functions $h_k$ may be treated analogously.

Under $P^{x_k}_{y_k}$, the process
$\{X_t, t\in[T_1, \zeta]\}$ 
is an $G_D(\, \cdot \, , y_k)$-process with the initial distribution
$$\mu_k(\, \cdot\,) \df
P^{x_k}_{y_k}(X(T_1)\in \,\cdot\,) =
P^{x_k} ( T_1 < T(D^c),
X(T_1) \in \,\cdot\,)G_D(\,\cdot\,,y_k)/
G_D(x_k,y_k),
$$
supported on $A_1$.
For a fixed $z\in A_1$, the process
$Y_t \df X_{\zeta -t}$ under $P^z_{y_k}$ has the
distribution $P^{y_k}_z$. 
If
$T_3 = \inf \{ t: Y_t \in A_2 \}$ then $T_3 = \zeta - T_2$.
The process $\{ Y_t , t\in [T_3, \zeta)\}$ under $P^z_{y_k}$
is a
$G_D(\, \cdot\, ,z)$-process with the initial distribution
$$\nu_k(\,\cdot\,) \df
P^{y_k} (T(A_2)< T(D^c), X(T(A_2)) \in \,\cdot\,)G_D(\,\cdot\,,z)/
G_D(y_k,z).
$$
For a fixed $v\in A_2$, the function
$y \to P^y(T(A_2)< T(D^c), X(T(A_2)) \in dv)$ is positive 
and harmonic in $D_2$
and vanishes on $ V_2 \cap \prt D$ and the same is true
for $z\to G_D(v ,z)$. By \Lbl{28},
$$
\frac{d\nu_k  }{d\nu_{3-k}} (v) =
\frac{
P^{y_k} (T(A_2)< T(D^c), X(T(A_2)) \in dv)G_D(v,z)G_D(y_{3-k},z)}
{G_D(y_k,z)
P^{y_{3-k}} (T(A_2)< T(D^c), X(T(A_2)) \in dv)G_D(v,z)}
\geq c_2.
$$
After reversing 
time again, 
we see that
the distributions of $X(T_2)$ under
$P^z_{y_1}$ and $P^z_{y_2}$
have 
Radon-Nikodym derivative bounded below by $c_2$.
The process $\{X_t, t\in[T_1, T_2]\}$ under
$P^z_{y_1}$ is a mixture of $h$-transforms converging
to $w$ with the mixing measure 
$P^z_{y_1}(X(T_2) \in dw)$ and the same remark applies to
$P^z_{y_2}$.
Hence, the distributions of $\{X_t, t\in[T_1, T_2]\}$ under
$P^z_{y_1}$ and $P^z_{y_2}$
have 
a 
Radon-Nikodym derivative bounded below by $c_2$.

We can prove in a similar way that
$d\mu_k (\,\cdot\,) /d\mu_{3-k}(\,\cdot\,) \geq c_1.$
The distributions of
$\{X_t, t\in[T_1, T_2]\}$ under
$P^{x_1}_{y_1}$ and $P^{x_2}_{y_2}$
have the Radon-Nikodym derivative bounded below by $c_1c_2$
because $P^{x_k}_{y_k}$ is a mixture of the measures
$P^z_{y_k}$ with the mixing measure $\mu_k$.\qed
\enddemo

\proclaim{\Tn{36} Lemma} Suppose that $f:\R^{d-1} \to \R$ is Lipschitz with
constant $\lambda$ and assume that $|f(x)| \leq 1$ for all $x$. Let
$$D = \{ x\in\R^d: |\wt x| < 1, f(\wt x) <x^d <2 \}.$$ 
There exists $c<\infty$ (which may depend on $\lambda$ but
does not otherwise depend on $f$) such that for
every $x\in \overline D$ and every positive harmonic function $h$ in $D$
$$
E^x_h \zeta < c.
\tag\tn{30}$$
\endproclaim
\demo{Proof} The result is essentially due to Cranston (1985) but
we refer the reader to the paper by Bass and Burdzy (1992).
Our domain $D$ is a special case of a ``twisted H\"older domain'' and
\Lbl{30} follows from Theorem 1.1 (i) (a) (C) of Bass and Burdzy (1992). A
direct inspection of its proof shows that $c$ depends only on the volume
and diameter of $D$ (under the assumption that $f$ is Lipschitz with
constant $\lambda$) and these quantities may be bounded 
independently of the particular form of $f$. \qed
\enddemo

\remark{\Tn{31} Remark} It is not necessary to assume in Lemma \lbl{36} that
$f$ is Lipschitz. It is enough to suppose that $f$ is upper
semicontinuous and $f(x)$ is bounded in the $L^p$-norm
for a suitable $p = p(d)$. 
This version of the result uses
Theorem 1.1 (i) (a) (A) of Bass and Burdzy (1992) which has a
considerably more complicated proof than Theorem 1.1 (i) (a) (C). We feel
it would not be fair to ask the reader to go through the former
proof in order to check that the constants may be chosen independently of
$f$.
\endremark

\proclaim{\Tn{32} Lemma} Suppose that $D\subset\R^d$ is a domain,
$x,y\in \overline D$, and for each $v=x,y$ there exist an 
orthonormal coordinate system
$CS_v$, a point $z_v\in D$, a Lipschitz function $f_v$ with constant 
$\lambda$ and a constant $c_v>0$ such that $|f_v| \leq c_v $,
$$\align
D_v \df  \{&z \in D: |\wt z| < c_v, -c_v < z^d < 2 c_v \ \ 
\text{in} \ CS_v \} \\
= \{&z \in \R^d: |\wt z| < c_v, f_v(\wt z) < z^d < 2 c_v \ \ 
\text{in} \ CS_v \} ,
\endalign$$
$$
z_v = (0,0, \dots , 0, 3c_v/2) \ \ \text{in}\ CS_v,
$$
$$
|\wt v| \leq c_v/2 \ \ \text{and} \ \ v^d\leq 3c_v/2 \ \ \text{in}\ CS_v,
$$
$$
D_x \cap D_y = \emptyset.
$$

If $E^{z_x}_{z_y} \zeta = c_1$ then
$$E^x_y \zeta \leq c_2 c_1 + c_3 (c_x^2 + c_y^2)$$
where $c_2$ and $c_3$ depend only on the dimension $d$ and
the Lipschitz constant $\lambda$.
\endproclaim
\demo{Proof} For $v=x,y$ let
$$\align
D^1_v & = \{ z\in D_v: |\wt z| < 3c_v/4, z^d < 7c_v/4 
\ \ \text{in}\ CS_v\}, \\
A_v & = \prt D^1_v \cap D, \\
Q_v &= \{ z \in \overline D_v : |\wt z | \leq c_v /2,
z^d \leq 3 c_v/2 \ \ \text{in}\ CS_v\}, \\
V_v &= \{ z\in \R^d : \dist (z, Q_v) < c_v / 8 \} .
\endalign$$
By the boundary Harnack principle 
\Lbl{81}(a), 
applied in $D_v$, assumption \Lbl{28} of
Lemma \lbl{27} holds. Let $T_1$ be the first hitting time of $A_x$ and $T_2$
be the last exit time from $A_y$. By Lemma \lbl{27},
$$
E^x_y (T_2 - T_1) \leq c_4 E^{z_x}_{z_y} (T_2 - T_1)
\leq c_4 E^{z_x}_{z_y} \zeta.
\tag\tn{33}$$
Lemma \lbl{36} and Brownian scaling \Lbl{80} imply that
$$E^x_y T_1 \leq c_5 c_x^2. \tag\tn{34}$$
The same lemma and time-reversal show that
$$E^x_y (\zeta -T_2) \leq c_5 c_y^2. \tag\tn{35}$$
The lemma follows from \Lbl{33}-\Lbl{35}. \qed
\enddemo

We now return to the specific domains, hypotheses, and notation 
of Theorem \slbl{1}{12}. 

\proclaim{\Tn{29} Lemma}  
Assume that $a<s_{j-1}<s_j<b $.
There exists $c_1>0$ such that for 
every positive harmonic function $h$ in $D_j$ which
vanishes on $\prt D_j \setminus \Lambda_{s_j}$ and every $x\in
\Lambda_{s_{j-1}}$,
$$E^x_h\zeta \geq c_1 f^2(s_{j-1}). $$

Moreover, there is a non-negative, non-constant and bounded
random variable $Y$ such that for every
$j$ and $x\in \Lambda_{s_{j-1}}$, the distribution of $\zeta$
under $P^x_h$ is stochastically larger than that of
$f^2(s_{j-1}) Y$.
\endproclaim
\demo{Proof} Let $B(y,r)$ denote the ball with center $y$ and radius $r$.
Let $c_2$ be the expected lifetime of conditioned Brownian motion in
$B(0,1)$ starting from 0 and converging to $x\in \prt B(0, 1)$. The
constant $c_2$ is strictly positive and does not depend on $x$ by
symmetry. For any harmonic function $g$ in $B(0,1)$, the $g$-process
starting from 0 is a mixture of processes conditioned to go to some
point of $\prt B(0,1)$ so its expected lifetime is also equal to
$c_2$. By scaling, the expected lifetime of any Brownian motion
conditioned by a harmonic function in $B(y,r)$ and starting from
$y$ is equal to $c_2 r^2$.

Let
$$\align
B_0 & = B((0,\dots,0,s_{j-1} + f(s_{j-1})/4),f(s_{j-1})/8), \\
T_1 & = \inf \{ t> T(B_0) : 
|X_t - X(T(B_0)) | 
= f(s_{j-1})/16 \}.
\endalign$$
Note that $B_0 \subset D_j$.
By the strong Markov property applied at $T(B_0)$,
$$
E^x_h \zeta \geq E^x_h [
(T_1 - T(B_0)) 
\jeden _{\{T(B_0) < \infty\} }]
= c_2 (f(s_{j-1})/16)^2 
P^x_h (T(B_0) < \infty ). 
\tag\tn{57}$$ 
Let $x_0 = (0,\dots,0,s_{j-1})$.
By Lemma
\lbl{26}, for all $x \in \Lambda_{s_{j-1}}$,
$$
P^x_h (T(B_0) < \infty) \geq c_3
P^{x_0}_{h} (T(B_0) < \infty). 
\tag\tn{16}$$
It is not hard to see that the constant $c_3$ may be chosen
independently of the particular form of $f$.
The probability 
$P^{x_0}_{h} (T(B_0) < \infty)$ 
is not less than
$$
P^{x_0} (T(B_0) < T(\prt D_j)) 
\inf_{y\in B_0} h(y) / h(x_0). 
$$
It is elementary to see that 
$P^{x_0} (T(B_0) < T(\prt D_j))$ 
is bounded below
and the usual Harnack principle shows that the same is
true for 
$\inf_{y\in B_0} h(y) / h(x_0)$. 
Hence, 
$P^{x_0}_{h} (T(B_0) < \infty)$ 
is bounded below by $c_4>0$ which together with 
\Lbl{57} 
and \Lbl{16} implies
$$E^x_h \zeta \geq c_2 (f(s_{j-1})/16)^2  c_3 c_4 .  $$

It is clear from our proof that $Y$ can be chosen as follows.
Let $\wt \zeta$ be the hitting time of $\prt B(0, 1/16)$ by a
Brownian motion starting from 0 
and let $W$ be an independent
random variable with $P(W=1)=1-P(W=0) = c_3c_4$. Then let
$Y=W Y'$, where $Y'= c_2 \min (\wt\zeta, 1)$.\qed 
\enddemo

\proclaim{\Tn{20} Lemma} 
Suppose that $ s_j < s_n$. Let $T_j^1 = T(\Lambda_{s_j})$ and 
$$\align
S_j^k &= \inf\{t> T_j^k : 
X_t \in \Lambda_{s_{j-1}} \cup \Lambda_{s_{j+1}} \},\quad k\geq 1,\\
T_j^k &= \inf\{t> S_j^{k-1} : 
X_t \in \Lambda_{s_j}  \},\quad k>1 .
\endalign$$
There exist $c_1<\infty $ and $p<1$
such that for all $k$ and for
every positive harmonic function $h$ in $D_n$ which
vanishes on $\prt D_n \setminus \Lambda_{s_n}$ and every $x\in D_n$
$$P^x_h (T_j^k < \infty) < c_1 p^k.$$
Moreover, if $i\geq 0$, $j+i<n$ and $x\in\Lambda_{s_{j+i}}$, then
$$
P^x_h (T_j^k < \infty) < c_1 p^{k+i}.
$$ 
\endproclaim
\demo{Proof}
Suppose $s_k<s_{k+1}\leq s_n$. We have
$$h(x) = \int_{\Lambda_{s_{k+1}} } h(y)
P^x(X(T(\Lambda_{s_{k+1}} )) \in dy)\tag\tn{22}$$
for $x \in \Lambda_{s_k}$.
The boundary Harnack principle implies that
$$\frac{P^{x_1}(X(T(\Lambda_{s_{k+1}} )) \in dy)}
{P^{x_2}(X(T(\Lambda_{s_{k+1}} )) \in dy)} 
\cdot\frac{P^{x_2}(T(\Lambda_{s_{k+1}} )<\infty)}
{P^{x_1}(T(\Lambda_{s_{k+1}} )<\infty)} 
<c_3<\infty\tag\tn{21}$$
for $x_1,x_2\in \Lambda_{s_k}$.
Let $z_k = (0,\dots,0,s_k)$.
It is easy to see that there is $c_4>0$ such that for all
$x \in \Lambda_{s_k}$ with $|\wt x| > (1-c_4) f(s_k)$, we have
$$P^x(T(\Lambda_{s_{k+1}} ) < \infty) 
< (c_3^{-1}/2) P^{z_k}(T(\Lambda_{s_{k+1}} ) < \infty).$$
This, \Lbl{22} and \Lbl{21} imply that 
$h(x) \leq h(z_k) /2$ for 
$x \in \Lambda_{s_k}$ with $|\wt x| > (1-c_4) f(s_k)$.
It follows that the maximum of $h$ on $\Lambda_{s_k}$ is attained at
a point in the set 
$$A_k \df \{x\in\Lambda_{s_k} : |\wt x| \leq (1-c_4) f(s_k)\}.$$

Let $a_k$ be the maximum of $h$ over $\Lambda_{s_k}$. Since
$$P^x(T(\Lambda_{s_{k+1}}) \leq T(\prt D_n)) < c_5<1$$ 
for $x \in \Lambda_{s_k}$, we have $a_k < c_5 a_{k+1}$ assuming
$a< s_k < s_{k+1} < b$.
It follows that $a_k < c_5^j a_{k+j}$. By the Harnack principle,
$h(x) > c_6 a_k$ for some $c_6>0$ and all $x\in A_k$. Let $m$ be
so large that $c_6 c_5 ^{-m} > 2$. Then $a_k < h(x)/2$
for all $x\in A_{k+m}$ provided $a< s_k < s_{k+m} < b$. We obtain
$$P^x_h(T(\Lambda_{s_j})<\infty) 
= \int_{\Lambda_{s_j}} \frac{h(y)}{h(x)} P^x(X(T(\Lambda_{s_j}))\in dy)
\leq 1/2\tag\tn{23}$$
for $x\in A_{j+m}$. Here and later in the proof we assume that
$a< s_j < s_{j+m} < b$. This assumption 
could 
be easily disposed of.
We have
$$P^{z_k} (T(A_{k+1}) < T(\prt D_n\cup \Lambda_{s_{k-1}})) > c_7>0$$
and an application of the Harnack principle shows that
$$P^{z_k}_h (T(A_{k+1}) < T( \Lambda_{s_{k-1}})) > c_8>0.$$
By Lemma \lbl{26}, 
$$P^x_h (T(A_{k+1}) < T( \Lambda_{s_{k-1}})) > c_9>0\tag\tn{24}$$
for all $x \in \Lambda_{s_k}$. By the strong Markov property applied
at the hitting times of $A_i$,
$$P^x_h (T(A_{j+m}) < T( \Lambda_{s_j})) > c_9^{m-1}\tag\tn{25}$$
for all $x \in \Lambda_{s_{j+1}}$. 
Let 
$$\align
U_1 &= \inf\{t>T(\Lambda_{s_{j+1}}): X_t \in A_{j+m}\},\\
U_2 &= \inf\{t>T(\Lambda_{s_{j+1}}): X_t \in \Lambda_{s_j}\},\\ 
U_3 &= \inf\{t> U_1 : X_t \in \Lambda_{s_j}\}.
\endalign$$
Then \Lbl{23}-\Lbl{25} imply that for $x\in \Lambda_{s_j}$
$$P^x_h (T_j^2 = \infty) 
\geq P^x_h(T(\Lambda_{s_{j+1}}) < T( \Lambda_{s_{j-1}}), 
U_1 < U_2, U_3 = \infty) > c_9^m/2 >0$$
for $x\in \Lambda_{s_j}$. 
Both conculsions of the lemma now follow 
by the repeated
application of the strong Markov property at the stopping times
$T_j^k$.\qed
\enddemo

\proclaim{\Tn{39} Lemma}
For all $x_1\in D_f$ such that 
$s_{k+1} \leq x_1^d \leq s_{k+2}$ and $x_2\in \Lambda_{s_k}$ we have
$E^{x_1}_{x_2} \zeta < c_1 f^2(s_k)$ where $E^{x_1}_{x_2}$ 
refers to the conditioned Brownian motion in $D_f$.
\endproclaim
\demo{Proof} We will suppose that $x_1 \in \Lambda_{s_{k+1}}$. The
modifications needed for the general case are obvious.

By Brownian scaling \Lbl{80}, we may assume that $f(s_k)=1$ and prove
that $E^{x_1}_{x_2} \zeta < c_1 $. Note that then $|x_1^d - x_2^d| =1/2$.

We have 
$$E^{x_1}_{x_2} \zeta = c_2 \int_{D_f} 
\frac{G_{D_f} ( x_1,z) G_{D_f} ( z,x_2)}{G_{D_f} ( x_1,x_2)} dz.$$

In view of Lemma \lbl{32} it will suffice to prove the lemma for
$x_1\in \Lambda_{s_{k+1}}$, $|\wt x_1| < c_3 f(s_{k+1})$, and
$x_2\in \Lambda_{s_k}$, $|\wt x_2| < c_3 $ for some $c_3 <1 $. 
Under
this additional assumption, $x_1$ and $x_2$ may be connected in
$D_f$ by a Harnack chain of balls of bounded length and this
implies that $G_{D_f} ( x_1,x_2) > c_4 >0$. Hence,
$$E^{x_1}_{x_2} \zeta < c_5 \int_{D_f} 
G_{D_f} ( x_1,z) G_{D_f} ( z,x_2) dz.\tag\tn{41}$$
Let
$$\align
A_j &= \{z\in D_f: |z- x_j| < 5,
|z - x_{3-j} | > |x_1-x_2| /2 \}, \quad 
j 
=1,2,\\
A_3 &= \{z\in D_f: |z- x_1| \geq 5, z^d < s_k \} ,\\
A_4 &= \{z\in D_f: |z- x_1| \geq 5, z^d > s_{k+1} \} .
\endalign$$

Assume for now that $d\geq 3$, and recall 
that $G(x,y) \df G_{\R^d} (x,y) = c_6 |x-y|^{2-d}$.
For $j=1,2$ we obtain 
$$\align
\int_{A_j} G_{D_f} ( x_1,z) G_{D_f} ( z,x_2) dz
&\leq \int_{A_j} G ( x_1,z) G ( z,x_2) dz\\
&\leq c_7 \int_{A_j} (|x_1-x_2| /2)^{2-d}
|z- x_j|^{2-d} dz  \tag\tn{40}\\
&\leq c_7 (|x_1-x_2| /2)^{2-d} \int_0^5 r^{2-d} r^{d-1} dr <
c_8<\infty.
\endalign$$

Let $x_0 = (0,\dots,0, s_k )$,
$$\align
\wt D &= \{x\in \R^d : x^d < s_k\},\\
D_* &= D_f \cup \{x\in \R^d : x^d \in (-\infty, s_k)
\cup (s_{k+1}, \infty) \},\\
M &= \{ x\in \wt D: |x-x_1| = 4\}.
\endalign$$
The Poisson kernel $K(x)$ in $\wt D$ with the pole at 
$x_0$ has the form 
$c_9 |x^d - s_k| /|x-x_0|^d$ (Doob (1984) 1.VIII.9). 
By the boundary Harnack principle,
$$G_{D_*}(x_1, x) \leq c_{10} K(x) $$
for $x\in M$ and, therefore, for all $x\in \wt D$
such that $|x-x_1| \geq 4$, in particular, for $x\in A_3$.
Hence, for $x\in A_3$,
$$G_{D_*}(x_1, x) \leq c_{11} |x^d - s_k| /|x-x_0|^d
\leq c_{11} |x-x_0|^{1-d}$$
and the same estimate holds for $G_{D_*}(x_2, x)$.
It follows that
$$\align
\int_{A_3} G_{D_f} ( x_1,z) G_{D_f} ( z,x_2) dz
&\leq \int_{A_3} G_{D_*} ( x_1,z) G_{D_*} ( z,x_2) dz\\
&\leq \int_{A_3} (c_{11} |z-x_0|^{1-d})^2 dz\tag\tn{42}\\
&\leq c_{12} \int _2^\infty r^{2(1-d)} r ^{d-1} dr < c_{13}
< \infty 
\endalign$$
and a similar estimate holds for $A_4$. Since 
$D_f \subset A_1 \cup A_2 \cup A_3 \cup A_4$, the lemma follows
from \Lbl{41}-\Lbl{42}. 

If $d=2$, an argument similar to the above could be given.
In this case, $\wt D$ should be replaced by a suitable wedge
with angle $\alpha<\pi$. The Green function in 
such a wedge decays
like $r^{-\pi/\alpha}$, and this is sufficient to make the 
bounding integrals finite. 
\qed
\enddemo

\proclaim{\Tn{48} Lemma} For $x\in D_f$ and $y\in\Lambda_{s_k}$, let 
$$
g_x^k (y) dy \df P^x_h (X(T(\Lambda_{s_k})) \in dy ).
$$
Then there exist $c_1 <
\infty$ and $c_2 < 1$ such that
$$
\frac
{g_{x_1}^n (y_1)}
{g_{x_1}^n (y_2)}
\geq a_i
\frac
{g_{x_2}^n (y_1)}
{g_{x_2}^n (y_2)}
\tag\tn{49}$$
and
$$a_i \geq 1 - c_1 c_2 ^i$$
for all $i>0$, all $n$, where
$x_1, x_2 \in D_{n-i}$ and 
$y_1, y_2 \in \Lambda_{s_n}$.
\endproclaim

\demo{Proof}
A standard application of the boundary Harnack
principle in the spirit of Lemma \lbl{26} shows that \Lbl{49}
holds for $i=1$ with some $a_1 >0$.

Assume that \Lbl{49} holds for all $n$ and for some $i$; we will show
that it holds for
$i+1$ as well. Let $j=n-i$.
By the strong Markov property applied at $T(\Lambda_{s_{n-1}})$,
$$
g_x^n (y)
= \int _{\Lambda_{s_{n-1}}}
g_x^{n-1} (v)
g_v^n (y)
dv
$$
for $y\in D_{j-1}$.
Now apply Lemma 6.1 of Burdzy, Toby and Williams (1989). Set in
that lemma
$V=W=\Lambda_{s_{n-1}}$ and $U=\emptyset$, set $f_1$ and $f_2$ equal
to our $g^{n-1}_{x_1}$ and $g^{n-1}_{x_2}$, set $g_z(v)$ equal to our
$g^n_v(z)$, and take $c=a_i$, $d=a_1$, and $b=1$.
The aforementioned
lemma implies that
$$
\frac
{g_{x_1}^n (y_1)}
{g_{x_1}^n (y_2)}
\geq a_{i+1}
\frac
{g_{x_2}^n (y_1)}
{g_{x_2}^n (y_2)}
$$
for all $y_1,y_2\in D_{j-1}$, where
$$a_{i+1} = a_i + a_1^2 (1- a_i).$$
Hence
$$1- a_{i+1} = 1 - a_i - a_1^2 (1- a_i)
= (1- a_i) (1-a_1^2)$$
and, by induction,
$$1- a_{i+1} \leq c_1 c_2 ^i ,$$
with  $c_2 \df 1-a_1^2<1$.\qed
\enddemo

\proclaim{\Tn{61} Corollary} With the notation of Lemma \lbl{48},
$$
a_{n-j}^{-1}\ge \frac{g^n_{x_1}(y)}{g^n_{x_2}(y)}\ge a_{n-j} 
$$
for every $j<n$, $x_1,x_2\in D_j$ and $y\in\Lambda_{s_n}$.
\endproclaim

\demo{Proof}
Let $M$ and $m$ be the supremum and infimum of $g^n_{x_1}(y)/g^n_{x_2}(y)$
over $y\in\Lambda_{s_n}$. 
By Lemma \lbl{48}, $m\ge a_{n-j}M$, and
$$
M g^n_{x_2}(y)\ge g^n_{x_1}(y)\ge m g^n_{x_2}(y).
$$
Integrating with respect to $y$ shows that $M\ge 1\ge m$, from 
which the desired conclusion follows. \qed
\enddemo 

\demo{Proof of Theorem \slbl{1}{12}} (i) 
We will first prove the lower bound in \SLbl{1}{17}. 

Suppose that $s_{j_0} \leq x^d < s_{j_0+1} < s_{j_0+2}<b$.
The other cases are left to the reader.
Let 
$T_j = T(\Lambda_{s_j}) $. For each $j>j_0+2$ the process
$\{X_t, t\in [T_{j-1}, T_j)\}$ under $P^x_h$
is a conditioned Brownian motion in $D_j$ starting from
a (random) point in $\Lambda_{s_{j-1}}$ and converging to
$\Lambda_{s_j}$ at its lifetime. 
By Lemma \lbl{29}, 
for $j\in[j_0+2,k_f -1]$,
$$ E^x_h(T_j - T_{j-1}) \geq c_1 f^2(s_{j-1})$$
and, therefore,
$$ E^x_h \zeta \geq \sum_{j=j_0 +2}^{k_f-1} E^x_h(T_j - T_{j-1})
\geq \sum_{j=j_0 +2}^{k_f-1} c_1 f^2(s_{j-1}) .\tag\tn{19}$$
Since 
$$c_2 f^2(s_{j-1}) < \int_{s_{j-1}}^{s_j} f(v) dv 
< c_3 f^2(s_{j-1}), $$
the sum on the right hand side of \Lbl{19} is bounded below by
$c_4 \int_{s_{j_0+1}}^{k_f-2}f(v) dv$.
Note that 
$$\int_{x^d}^{s_{j_0+1}}f(v) dv< c_5 
\int_{s_{j_0+1}}^{s_{j_0+2}}f(v) dv$$
and
$$\int_{k_f-2}^b f(v) dv < c_5 \int_{k_f-3}^{k_f-2}f(v) dv .$$
Hence
$$\int_{x^d}^b f(v) dv < c_6 \int_{s_{j_0+1}}^{k_f-2}f(v) dv$$
and, therefore,
$$E^x_h \zeta \geq c_7 \int_{x^d}^b f(v) dv .$$

(ii) Next we will prove (ii) of Theorem \slbl{1}{12}. 

First note that $k_f=\infty$. 
Recall the definitions of $j_0$ and the $T_j$'s from part (i) of the proof.
By Lemma \lbl{29} and the strong Markov property applied at 
$T_j$'s, there exist non-negative (not necessarily independent) random
variables $Z_j$ and i.i.d. non-negative random variables $Y_j$
such that 
$$\sum_{j=j_0 +2}^{\infty} (T_j - T_{j-1})\tag\tn{46}$$ 
has the same distribution as 
$$\sum_{j=j_0 +2}^{\infty} (Z_j + f^2(s_{j-1}) Y_j).\tag\tn{47}$$
For later use, note that, as in the proof of Lemma \lbl{29},
we can write  $Y_j=W_j Y_j'$, where
the $Y_j'$ are independent of the $Z$'s and $W$'s, with some common mean
$\mu$ and variance $\sigma^2$. Each $W_j$ takes values $0$ or $1$, and
$W_j=1$ with some common probability $p$, even if conditioned on the
preceding $W$'s and on $\{X_t, t\in[0, T_{j-1}]\}$. Thus the $W_j$ are
i.i.d., though they may not be independent of the $Z_j$.

It is elementary to check that 
$\sum_{j=j_0 +2}^{\infty}  f^2(s_{j-1}) =\infty$
because $\int_{x^d}^b f(v) dv =\infty$. Hence,
$$\sum_{j=j_0 +2}^{\infty}  E (f^2(s_{j-1})Y_j) =\infty.$$
Recalling that each $Y_j$ is
non-negative, non-constant and bounded,
the three-series theorem 
now easily implies that a.s.
$$\sum_{j=j_0 +2}^{\infty} f^2(s_{j-1})Y_j =\infty.$$
It follows that the sums in \Lbl{47}, and therefore in \Lbl{46},
must be infinite a.s. 

(iii) We are going to prove the lower bound in \SLbl{1}{18}. 

Let $j_0$, the $Y_j$'s, etc.  be as in part (ii) of the proof. 
By adjusting the first and last $Z$, if necessary, we can guarantee
that 
$$
\align
\zeta&=\sum_{j=j_0 +2}^{k_f-1}(Z_j+f^2(s_{j-1}) Y_j)\\
&=\sum_{j=j_0 +2}^{k_f-1}(Z_j+f^2(s_{j-1}) \mu W_j) +
\sum_{j=j_0 +2}^{k_f-1}(f^2(s_{j-1})W_j(Y_j'-\mu)).\tag\tn{77}
\endalign
$$
Therefore by independence,
$$
\align
\var_h^x \zeta &=
\var_h^x\left(\sum_{j=j_0 +2}^{k_f-1}(Z_j+f^2(s_{j-1}) \mu W_j)\right) 
+ \sum_{j=j_0 +2}^{k_f-1}E_h^x((f^2(s_{j-1})W_j(Y_j'-\mu))^2)\\
&\geq \sum_{j=j_0 +2}^{k_f-1}E_h^x((f^2(s_{j-1})W_j(Y_j'-\mu))^2)\\
&\ge \sum_{j=j_0 +2}^{k_f-1} f^4(s_{j-1})p\sigma^2
\ge c_3\int_{x^d}^b f^3(v) dv.
\endalign 
$$ 

(iv) We will now prove part (v) of Theorem \slbl{1}{12}. 

We will again invoke the
$Y_j$'s and $Z_j$'s of part (ii) of the proof.
Suppose that $\int_{x^d}^b f^3(v) dv = \infty$. Then necessarily
$b=\infty$. Let us assume that 
$$\limsup_{v\to\infty} f(v) <\infty.\tag\tn{64}$$
In order to simplify the notation, suppose that $x^d = s_{j_0}$.

First, let $w_1$, $w_2$, \dots be any sequence of 0's and 1's, 
such that 
$$
\sum_{j>j_0}f^4(s_{j-1})w_j=\infty.
$$ 
Consider
$$
\wt Y_k =  \sum_{j=j_0 +1}^k  f^2(s_{j-1}) w_j(Y_j'-\mu)
\quad\text{and}\quad\wh Y_k = \wt Y_k/ (\var \wt Y_k)^{1/2}.
$$
Since the $Y_j'$s are uniformly bounded,
the Lindeberg-Feller condition can be easily verified using
\Lbl{64} and it follows that the distributions of
$\wh Y_k$ converge to the standard normal distribution as 
$k\to\infty$. In fact it is simple to show, using \Lbl{64} and the
Berry-Eseen theorem, that for every $c_1<\infty$ and $c_2>0$
there exists a $c_3<\infty$ such that 
$$
\text{$P(\wt Y_k \in (u, u + c_1) ) < c_2/2$\ \ 
for every $u\in\R$, if \ \ $\var\wt Y_k >c_3$.}\tag\tn{65}
$$

Since $\sum_{j>j_0} f^4(s_{j-1})W_j=\infty$ almost surely, we can
choose a $k_0<\infty$ such that 
$$
P^x_h\left(\sum_{j=j_0+1}^kf^4(s_{j-1})W_j>c_3\right)>1-c_2/2
$$
for every $k\ge k_0$. Also, as in \Lbl{77} we have that
$$
T(\Lambda_{s_k})
=\sum_{j=j_0 +2}^{k}(Z_j+f^2(s_{j-1}) \mu W_j) +
\sum_{j=j_0 +2}^{k}(f^2(s_{j-1})W_j(Y_j'-\mu)).
$$
Therefore, conditioning on the values
of $W_j, j>j_0$ yields  that
$$
P^x_h(T(\Lambda_{s_k})\in(u,u+c_1))<c_2
$$
for every $u\in\R$. 

The case when \Lbl{64} fails is not hard and is left to the reader.

(v) Next we prove the upper bound in \SLbl{1}{17}. 

Suppose that $s_{n+1} \leq x^d \leq s_{n+2}$.
Let $L$ be the last exit time from $\Lambda_{s_n}$. Under $P^x_h$,
the process
$\{X_t, t \in [0,L]\}$ is a conditioned Brownian motion in $D_f$
starting from $x$ and converging to a (random) point of 
$\Lambda_{s_n}$. Lemma \lbl{39} implies that 
$E^x_h L <  c_1 f^2(s_n)$ and this in turn implies that
$$E^x_h L < c_2 \int _{x^d}^{s_{n+3}} f(v) dv.\tag\tn{43}$$

For every $\eps>0$, the process $\{X_{t+L+\eps}, t\ge 0\}$ 
under $P^x_h$ is an
$h$-process in the domain $D_g$ where 
$g(s) = f(s) \jeden_{(s_n,\infty)}(s)$. This and \Lbl{43} show
that 
\SLbl{1}{17} 
will follow once we prove that
$$E^x_h \zeta < c_3 \int_a^b f(v) dv.$$

Let $M_k = \{ y\in D_f: s_{k-1} < y^d < s_{k+1} \} $
and consider an $h_0$-process in $M_k$ for some positive
harmonic function $h_0$ in $M_k$. A variation of
Lemma \lbl{36} shows that 
$$E^y_{h_0} \zeta < c_4\tag\tn{50}$$
for all
$y\in M_k$, 
provided $f(s_k) =1$. By scaling, 
$$E^y_{h_0} \zeta < c_4 f^2(s_k)\tag\tn{54}$$
for any value of $f(s_k)$.

Recall 
the 
stopping times $S_j^k$ and $T_j^k$ from Lemma \lbl{20}
and let $F_j^k \df \{T_j^k < \infty\}$.
Let $T_0$ be the hitting time of $\bigcup_k \Lambda_{s_k}$.
We have 
$$\zeta = T_0 + \sum_{j,k}  (S_j^k - T_j^k)\jeden_{F_j^k}.
\tag\tn{44}$$
Given $T_j^k < \infty$, the process
$\{X_t, t \in [T_j^k, S_j^k] \}$ is a conditioned Brownian
motion in $M_k$ and, therefore, 
$$E^x_h [ (S_j^k - T_j^k) \mid F_j^k ] < c_4 f^2(s_j).$$ 
By Lemma \lbl{20},
$$\sum_k E^x_h (S_j^k - T_j^k)\jeden_{F_j^k}
< c_5 f^2(s_j).\tag\tn{45}$$
Recall that $s_{n+1} \leq x^d \leq s_{n+2}$. Hence
$E^x_h T_0 < c_4 f^2 (s_n)$. This and \Lbl{44}-\Lbl{45} yield
$$E^x_h \zeta \leq c_6 \sum _j f^2(s_j).$$
It is easy to check that the last quantity is bounded by
$c_7 \int_a^b f(v) dv$.

(vi) We will 
now prove the upper bound for the variance in \SLbl{1}{18}. 
Recall $M_k$ and 
the use of 
an $h_0$-process in $M_k$
from part (v) of the proof. The Chebyshev
inequality and \Lbl{50} show that
$P^x_{h_0} (\zeta > c_1 ) < c_2$ for some $c_1<\infty$, $c_2 < 1$
and all $x\in M_k$ provided $f(s_k) = 1$. By the Markov property
applied repeatedly at the multiples of $c_1$,
$P^x_{h_0}(\zeta > jc_1) < c_2^j$.
Hence
$E^x_{h_0} \zeta^2 < c_3 $
in the case $f(s_k) = 1$ and, by scaling,
$$E^x_{h_0} \zeta^2 < c_3 f^4(s_k)\tag\tn{55}$$
for any value of $f(s_k)$, all $x\in M_k$ and all harmonic
functions $h_0$ in $M_k$.

Let $S_j^k$ and $T_j^k$ be as in Lemma \lbl{20}. 
Let $F_j^k \df \{T_j^k < \infty\}$.
Given 
$F_j^k$, the process
$\{X_t, t \in [T_j^k, S_j^k] \}$ is a conditioned Brownian
motion in $M_k$ and this implies in view of \Lbl{54} and
\Lbl{55}, that 
$$
\aligned
E^x_h [(S_j^k - T_j^k) \mid F_j^k]&< c_4 f^2(s_j)\quad\text{and}\\
E^x_h [(S_j^k - T_j^k)^2 \mid F_j^k ] &< c_3 f^4(s_j).
\endaligned
\tag\tn{53}
$$ 

Let $\Theta_j^k \df (S_j^k - T_j^k)\jeden_{F_j^k}$. 
Define $q$ by the condition that $s_{q-1}<x^d\le s_q$, and recall
from Lemma \lbl{20} that
$$
P^x_h(F_j^k) \le 
\cases
c_5 c_6 ^{k+q-j},& j<q\\
c_5 c_6 ^k,& j\ge q,
\endcases
\tag\tn{59}
$$
where $c_6<1$. This and \Lbl{53} imply that
$$
E^x_h [\Theta _j^k] \le 
\cases
c_4 c_5 c_6^{k+q-j}f^2(s_j),& j<q\\
c_4 c_5 c_6^k f^2(s_j),& j\ge q,
\endcases
\tag\tn{58}
$$
and
$$
E^x_h [(\Theta _j^k)^2] \le 
\cases
c_3 c_5 c_6^{k+q-j}f^4(s_j),& j<q\\
c_3 c_5 c_6^k f^4(s_j),& j\ge q,
\endcases
\tag\tn{60} 
$$

Now assume that
$j<n$, and let 
$$
A=\{T^k_j<T^1_n\},\quad B=\{T^1_n<T^k_j\}, 
\quad B_i=\{T^{i-1}_j<T^1_n<T^i_j\}
$$
where $T^0_j$ is taken to be $0$. Then
$$
\aligned
\cov^x_h(\Theta^k_j,\Theta^m_n) 
&=E^x_h((\Theta^k_j-E^x_h\Theta^k_j)
(\Theta^m_n-E^x_h\Theta^m_n))\\
&=E^x_h((\Theta^k_j-E^x_h\Theta^k_j)
(\Theta^m_n-E^x_h\Theta^m_n)\jeden_A) +\\
&\qquad +
E^x_h((\Theta^k_j-E^x_h\Theta^k_j)
(\Theta^m_n-E^x_h\Theta^m_n)\jeden_B)\\
&\df \text{\sl I} + \text{\sl II}.
\endaligned
\tag\tn{51}
$$

Consider term $\text{\sl I}$ of \Lbl{51}. If $q>n$ then $\text{\sl
I}=0$ automatically. So suppose that $q\le j$. By Corollary \lbl{61}
and the strong Markov property at $T^1_n$,
$$
|E^y_h\Theta^m_n - E^x_h\Theta^m_n|\le 
c_7 c_8^{n-j} E^x_h\Theta^m_n
$$
for any $y\in D_j$, where $c_8<1$. In particular, 
$$
|E^x_h(\Theta^m_n\mid \Cal F_{S_j^k})
- E^x_h\Theta^m_n | \le c_7 c_8^{n-j} E^x_h(\Theta^m_n)
$$
on $A$. Thus, by \Lbl{58},
$$\align
\text{\sl I}&=E^x_h\left[(\Theta_j^k - E^x_h\Theta_j^k )\jeden_A
E^x_h(\Theta^m_n - E^x_h(\Theta^m_n)\mid\Cal F_{S_j^k}) \right]\\
&\leq
E^x_h\left[ |\Theta _j^k - E^x_h\Theta _j^k |\cdot \jeden_A \cdot
|E^x_h(\Theta^m_n\mid \Cal F_{S_j^k})
- E^x_h\Theta^m_n | \right]\\
&\le 2 c_7 c_8^{n-j} E^x_h(\Theta^m_n) E^x_h(\Theta^k_j)
\le c_9 c_8^{n-j} c_6^{k+m} f^2(s_j)f^2(s_n).
\endalign
$$

If, on the other hand, we have $j<q\le n$, then by a similar argument,
$$
|E^x_h(\Theta^m_n\mid \Cal F_{S_j^k})
- E^x_h\Theta^m_n | \le c_7 c_8^{n-q} E^x_h(\Theta^m_n)
$$
on $A$, and 
$$
\text{\sl I}\le 2 c_7 c_8^{n-q} E^x_h(\Theta^m_n) E^x_h(\Theta^k_j)
\le c_9 c_8^{n-q} c_6^{k+m+q-j} f^2(s_j)f^2(s_n).
$$
Taking $c_{11}=\max(c_8,c_6)$, it follows that 
$$
\text{\sl I}\le c_9 c_{11}^{n-j} c_6^{k+m}
f^2(s_j)f^2(s_n),\tag\tn{56}
$$
regardless of the value of $q$.

Consider now the term {\sl II} of \Lbl{51}.
By \Lbl{60}, and by Lemma \lbl{20} again, 
$$
\align
E^x_h((\Theta^k_j)^2\jeden_B) 
&= \sum_{i=1}^k E^x_h((\Theta^k_j)^2\jeden_{B_i}) \\
&= \sum_{i=1}^k E^x_h(E^x_h((\Theta^k_j)^2\jeden_{B_i}\mid\Cal
F_{S^1_n}))\\
&\le \sum_{i=1}^k f^4(s_j) c_0 c_3 c_6^{n-j+k-i+1}P^x_h(B_i)\\
&\le \sum_{i=1}^k f^4(s_j) c_0
c_3 c_6^{n-j+k-i+1}P^x_h(F^{i-1}_j)\\   
&\le k c_0^2 c_3 c_6^{n-j+k} f^4(s_j) 
\le c_{12} c_{13}^{n-j+k} f^4(s_j),
\endalign
$$
where $c_{13}<1$. As a result, 
$$
\align
\text{\sl II}
&\le (E^x_h((\Theta^k_j)^2\jeden_B))^{1/2}
(E^x_h((\Theta^m_n)^2))^{1/2}\\
&\le f^2(s_j)f^2(s_n) (c_{12} c_{13}^{n-j+k} c_0 c_3 c_6^m)^{1/2}\\
&\le c_{14} c_{15}^{k+m+n-j} f^2(s_j)f^2(s_n),
\endalign
$$
where $c_{14}<1$.
Combining this with \Lbl{51} and \Lbl{56}, it follows that
$$
\cov^x_h(\Theta^k_j,\Theta^m_n)\le c_{16} c_{17}^{k+m+|n-j|} 
f^2(s_j)f^2(s_n),\tag\tn{52}
$$
for $j<n$, where $c_{17}<1$. By symmetry, the same is true for $j>n$,
and the inequality is even simpler to prove if $j=n$ (\Lbl{56} is no 
longer needed). Thus, \Lbl{52} holds for every $j,k,m,n$.

If $\int_a^b f^3(v) dv = \infty$ then
the upper bound in \SLbl{1}{18} is trivial. Assume therefore that 
$\int_a^b f^3(v) dv < \infty$. Then for each $\eps>0$ there 
are only finitely many
$j$ such that $f(s_j) > \eps$. Hence we may
choose an ordering $\{j_i\}_{i\geq1}$ of the set
$\{k: a<s_k<b\}$ which satisfies
$f(s_{j_{i+1}}) \leq f(s_{j_i})$ for all $i$.
By \Lbl{52}
$$
\aligned
\var^x_h \zeta &= \var^x_h 
\biggl( \sum_{j,k} \Theta _j^k \biggr)
= \sum_{j,k,n,m} \cov^x_h (\Theta _j^k, \Theta _n^m)\\
&\le 2 \sum _i \sum _{n\ge i} \sum _k \sum _m
\cov^x_h (\Theta _{j_i}^k, \Theta _{j_n}^m)\\
&\le 2 \sum _i \sum _{n\ge i} \sum _k \sum _m
c_{13} c_{12}^{k+m+|j_n-j_i|}  f^2(s_{j_i}) f^2(s_{j_n})\\ 
&\le \sum _i \sum _{n\ge i} c_{14} c_{12} ^{|j_n-j_i|}f^4(s_{j_i})\\
&\leq \sum_j c_{15} f^4(s_j)
\leq c_{16} \int_a^b f^3(v) dv.
\endaligned \tag\tn{67}
$$ 

(vii) Next we will prove part (iv) of Theorem \slbl{1}{12}.

Fix some $x\in D_f$ and
suppose for convenience that $x^d = s_q$ for some $q$.
Recall $S_j^k, T_j^k, F_j^k$ and $\Theta _j^k$ from part (v) of
the proof. With slightly more work, the argument for \Lbl{52}
can be seen to yield the following improved estimate:
$$
\cov^x_h(\Theta^k_j,\Theta^m_n)\le
\cases
c_1 c_2^{k+m+|n-j|} f^2(s_j) f^2(s_n), &j,n\ge q\\
c_1 c_2^{k+m+|n-j|} c_3^{q-j} f^2(s_j) f^2(s_n), &j< q\le n\\
c_1 c_2^{k+m+|n-j|} c_3^{q-j}c_3^{q-n} f^2(s_j) f^2(s_n), &j,n< q,
\endcases
$$
where $c_2, c_3<1$.

Now we assume that the Lipschitz constant of $f$ is so small
that for each $j$,
$$
\frac{f^2(s_{j-1})}{f^2(s_j)} < 
\frac{c_3^{-1}+1}{2}.
$$
Therefore 
$$
\cov^x_h(\Theta^k_j,\Theta^m_n)\le
\cases
c_1 c_2^{k+m+|n-j|} f^2(s_j) f^2(s_n), &j,n\ge q\\
c_1 c_2^{k+m+|n-j|} c_4^{q-j} f^2(s_q) f^2(s_n), &j< q\le n\\
c_1 c_2^{k+m+|n-j|} c_4^{q-j} c_4^{q-n} f^4(s_q), &j,n<q,
\endcases\tag\tn{62}
$$
for some $c_4<1$. 

If $\int_{x^d}^b f^3(v) dv = \infty$ then
\SLbl{1}{37} obviously holds. Assume that 
$\int_{x^d}^b f^3(v) dv < \infty$. Then we may
choose an ordering $\{j_i\}_{i\geq1}$ of the set
$\{k: x^d \leq s_k<b\}$ which satisfies
$f(s_{j_{i+1}}) \leq f(s_{j_i})$ for all $i$. 
Let $j_{i_0}=q$. Then
in view of \Lbl{62},
$$
\aligned
\var^x_h \zeta &= \var^x_h 
\biggl( \sum_{j,k} \Theta _j^k \biggr)
= \sum_{j,k,n,m} \cov^x_h (\Theta _j^k, \Theta _n^m)\\
&\le 2 \sum _i \sum _{n\ge i} \sum _k \sum _m
\cov^x_h (\Theta _{j_i}^k, \Theta _{j_n}^m)\\
&\qquad +2\sum_{j\le n<q}\sum_k\sum_m\cov^x_h 
(\Theta _j^k, \Theta_n^m)\\ 
&\qquad +2\sum_{j<q}\sum_{i\ge i_0}\sum_k\sum_m
\cov^x_h (\Theta _j^k,\Theta_{j_i}^m)\\
&\qquad +2\sum_{j<q}\sum_{i< i_0}\sum_k\sum_m
\cov^x_h (\Theta _j^k,\Theta_{j_i}^m)\\
&\le 2 \sum _i \sum _{n\ge i} \sum _k \sum _m
c_1 c_2^{k+m+|j_n-j_i|}  f^2(s_{j_i}) f^2(s_{j_n})\\ 
&\qquad +2\sum_{j\le n<q}\sum_k\sum_m
c_1 c_2^{k+m+|n-j|} c_4^{q-j} c_4^{q-n} f^4(s_q)\\
&\qquad +2\sum_{j<q}\sum_{i\ge i_0}\sum_k\sum_m
c_1 c_2^{k+m+|j_i-j|} c_4^{q-j} f^2(s_q) f^2(s_{j_i})\\
&\qquad +2\sum_{j<q}\sum_{i< i_0}\sum_k\sum_m
c_1 c_2^{k+m+|j_i-j|} c_4^{q-j} f^2(s_q) f^2(s_{j_i})\\
&\le \sum_i c_5 f^4(s_{j_i}) +c_6 f^4(s_q) +c_7 f^4(s_q)
+\sum_i c_8 f^4(s_{j_i})\\
&\le c_9\sum_{j\ge q} f^4(s_j)
\le c_{10} \int_{x^d}^b f^3(v)dv. \qed
\endaligned\tag\tn{63}
$$
\enddemo 

Because they use similar arguments to those just given, we include
the following two subsidiary results in this section. 

\proclaim{\Tn{66} Corollary} Suppose that $D_f$ and $h$ are
as in Theorem \slbl{1}{12}. 
Assume that \break
$\int_a^b f^3(v) dv < \infty$. Then
$$\lim_{x^d \to\infty}
\sup \{ \var^x_h T(\Lambda_{u}) : u > x^d\} = 0.\tag\tn{69}$$
\endproclaim
\demo{Proof} Recall 
the notation from the proof of Theorem \slbl{1}{12}. 
As in the proof of \Lbl{67}, 
for every $x$ and for every $u=s_i$,
$$\var^x_h T(\Lambda_{u}) =
\sum_{j,k,n,m}\cov^x_h (\Theta _j^k\jeden_{\{T^k_j<T(\Lambda_u)\}},
\Theta _n^m\jeden_{\{T^m_n<T(\Lambda_u)\}}).
$$
An examination of the proof of \Lbl{67} shows that
the terms of this sum are bounded by the terms of an absolutely
convergent series, uniformly in $x$ and in $u=s_i$. With a little
more work, it is easy to see that this domination holds 
for $u\in(a,b)$ as well. 
For fixed $j,k,m$ and $n$,
$$
\cov^x_h (\Theta _j^k\jeden_{\{T^k_j<T(\Lambda_u)\}},
\Theta _n^m\jeden_{\{T^m_n<T(\Lambda_u)\}})\to 0
$$ 
as $x^d\to \infty$,
uniformly in $u$,
because of \Lbl{59}. This easily implies \Lbl{69}.\qed 
\enddemo

\proclaim{\Tn{116} Lemma} 
Assume that $D_f$ and $h$ are as in Theorem \slbl{1}{12}. Set 
$$
f_*(v) \df \sup_{u\geq v} f(u).
$$ 
There exists a $c_1<\infty$ 
such that for all $u$ and all
$x_1,x_2\in D_f$ with $x_1^d = x_2^d < u$ we have
$$|E^{x_1}_h T(\Lambda_u) - E^{x_2}_h T(\Lambda_u)|
\leq c_1 f^2_*(x_1^d).$$
\endproclaim
\demo{Proof}
We will use an argument from part (v) of the proof of Theorem
\slbl{1}{12}. Suppose that $s_{n+1} \leq x_1^d \leq s_{n+2}$ and let
$L$ be the last exit from $\Lambda_{s_n}$. It has been proved that
$$E^{x_k}_h L < c_2 f^2 (s_n)\tag\tn{117}$$
for $k=1,2$ (see the paragraph 
preceding \Lbl{43}). 
Recall the definitions of $T_0, S_j^k, T_j^k$ and $F_j^k$ from the
same proof, and set 
$$
G^k_j\df F^k_j\cap\{T^k_j<T(\Lambda_u)\}. 
$$ 
We have
$$E^{x_k}_h T_0 < c_3 f^2 (s_{n+1})\tag\tn{119}$$
by an argument analogous to that proving 
\Lbl{54}. By Lemma \lbl{48}, 
the Radon-Nikodym derivative of the
initial distributions of
$\{ X_t, t\in [T_j^k , S_j^k] \}$ under 
$P^{x_1}_h(\,\cdot\mid~G^k_j)$ and $P^{x_2}_h(\,\cdot\mid G^k_j)$ 
differs from 1 by no more than
$c_4 c_5^{|n-j|}$ where $c_5 < 1$. It follows that
$$
|E^{x_1}_h [(S_j^k - T_j^k) \jeden_{G_j^k}]
- E^{x_2}_h [(S_j^k - T_j^k) \jeden_{G_j^k}]|
\leq c_4 c_5^{|n-j|} E^{x_1}_h [(S_j^k - T_j^k) \jeden_{G_j^k}].$$ 
Now \Lbl{45} implies that
$$
\aligned
\bigg |\sum_k E^{x_1}_h (S_j^k - T_j^k) &\jeden_{G_j^k}
- \sum_k E^{x_2}_h (S_j^k - T_j^k) \jeden_{G_j^k}\bigg|\\
&\leq c_4 c_5^{|n-j|} \sum_k 
E^{x_1}_h (S_j^k - T_j^k) \jeden_{G_j^k}\\
&\leq c_4 c_5^{|n-j|} c_6 f^2(s_j).
\endaligned \tag\tn{118}
$$ 
Since 
$$
\sum_{k\geq 1}(S_j^k - T_j^k) \jeden_{G_j^k}\le
T(\Lambda_u) \leq T_0 + L + \sum_{j\geq n} \sum_{k\geq 1}
(S_j^k - T_j^k) \jeden_{G_j^k},
$$ 
we obtain from \Lbl{117}-\Lbl{118} that 

$$
\align
|E^{x_1}_h T(&\Lambda_u) - E^{x_2}_h T(\Lambda_u)|\\
&\leq 2 c_2 f^2 (s_n) + 2 c_3 f^2 (s_{n+1}) +
\sum_{j\geq n} c_4 c_5^{|n-j|} c_6 f^2(s_j)\\ 
& \leq 2 c_2 f_*^2 (s_n) + 2 c_3 f_*^2 (s_n) +
\sum_{j\geq n} c_4 c_5^{|n-j|} c_6 f_*^2(s_n) 
\leq c_7 f_*^2(s_n) .\qed
\endalign
$$ 
\enddemo


\sectno = 3
\tagno=1

\subheading{3. Disintegration of harmonic functions}

The purpose of this section is to prove Theorem \SLbl{1}{14}. 
Unless otherwise indicated, the notation and general hypotheses of
Theorem \SLbl{1}{14} will be assumed throughout this section. 

Fix some $x_0\in D_f$ and let
$g(u) \df - E^{x_0}_h T(\Lambda_u)$. Recall that
$f_*(v) = \sup_{u\geq v} f(u)$.
Note that in either case (a) or (b) of Theorem \SLbl{1}{14} (i),
we have that $f(v) \to 0$ as $v\to\infty$. 

\proclaim{\Tn{148} Lemma} 
Suppose that one of the assumptions (a) or (b) of Theorem \SLbl{1}{14}
(i) is satisfied. Then 
$$
\lim_{u\to\infty} (T(\Lambda_u)+g(u))\quad
\text{exists } P^{x_0}_h\text{-a.s.}
$$
\endproclaim
\demo{Proof} 
Lemma \slbl{2}{116} and Corollary \SLbl{2}{66} show that
for $k\geq 1$, we can choose $u_k$ such that
$$|E^{x_1}_h T(\Lambda_u) - E^{x_2}_h T(\Lambda_u)|
\leq c_1 f^2_*(u_k) \leq 1/k^2\tag\tn{120}$$
for all $x_1,x_2\in D_f$ and $u$ with $u_k \leq x_1^d=x_2^d <
u$. We may also assume that
$$\var^x_h T(\Lambda_u) \leq 1/k^6\tag\tn{121}$$ for $x\in D_f$
and $u$ with $u_k \leq x^d < u$. 

Suppose $u\in [u_k, u_{k+1})$. Since 
$$T(\Lambda_{u_{k+1}}) = (T(\Lambda_{u_{k+1}}) - T(\Lambda_u))
+ T(\Lambda_u),$$
we have
$$ g(u_{k+1}) = - E^{x_0}_h(T(\Lambda_{u_{k+1}}) - T(\Lambda_u))
+g(u).$$
This, \Lbl{120}, and the strong Markov property applied at
$T(\Lambda_u)$ imply that
$$|E^x_hT(\Lambda_{u_{k+1}}) 
+(g(u_{k+1}) - g(u)) | \leq 1/k^2\tag\tn{122}$$
for all $x\in D_f$ such that $x^d = u$.
The Chebyshev inequality and \Lbl{121} yield that
$$P^x_h(|T(\Lambda_{u_{k+1}})  - E^x_hT(\Lambda_{u_{k+1}}) |
\geq 1/k^2) \leq k^4 \var^x_h T(\Lambda_{u_{k+1}}) \leq 1/k^2,$$
if $x^d = u$. This and \Lbl{122} give
$$P^x_h(|T(\Lambda_{u_{k+1}})  +(g(u_{k+1}) - g(u)) | 
\geq 2/k^2) \leq 1/k^2,$$
for $x\in D_f$ such that $x^d = u$. By the strong Markov
property applied at $T(\Lambda_u)$, 
$$
P^x_h 
(|T(\Lambda_{u_{k+1}})- T(\Lambda_u)
+(g(u_{k+1}) - g(u)) | 
\geq 2/k^2) \leq 1/k^2,\tag\tn{125}
$$
for any $x\in D_f$ with $x^d\le u$. In particular, 
$$P^x_h(|T(\Lambda_{u_{k+1}})- T(\Lambda_{u_k}) 
+(g(u_{k+1}) - g(u_k)) | 
\geq 2/k^2) \leq 1/k^2,\tag\tn{126}$$
if $x^d\le u_k$. 

Fix some $c_2 > 0$ and find $j_0$ so large that
$\sum_{j\geq j_0} 2/j^2 < c_2$. Suppose that 
$k>j_0$, $x^d\le u$, 
and recall
that $u\in [u_k, u_{k+1})$. Then \Lbl{125}-\Lbl{126} imply that with
$P^x_h$-probability 
larger than 
$1-c_2$, the event
$$
\multline
\{|T(\Lambda_{u_{k+1}})- T(\Lambda_u)
+(g(u_{k+1}) - g(u)) | \leq 2/k^2\}\\
\cap\bigcap_{j\ge j_0}\{|T(\Lambda_{u_{j+1}})- T(\Lambda_{u_j})
+(g(u_{j+1}) - g(u_j)) | \leq 2/j^2\}
\endmultline\tag\tn{123}
$$
occurs. Let 
$$A_v \df \{ |(T(\Lambda_{u_m}) + g(u_m)) - (T(\Lambda_v) +g(v))| 
< c_2 \ \ \forall u_m \geq v \}.$$
If the event in \Lbl{123} holds 
then $A_u$ holds,
because in such a case we have
$$\align
|(T(\Lambda_{u_m}) + g(u_m)) &- (T(\Lambda_u) +g(u))|\\
& \leq |(T(\Lambda_{u_{k+1}}) + g(u_{k+1})) - (T(\Lambda_u) +g(u))| \\
&\qquad + \sum_{j=k+1}^{m-1}
|(T(\Lambda_{u_{j+1}})- T(\Lambda_{u_j}))
+  (g(u_{j+1}) - g(u_j))| \\
&\leq 2/ k^2 +\sum_{j=k+1}^{m-1} 2/j^2 < c_2.
\endalign$$
Hence $P^x_h(A_u) > 1-c_2$. 

Let 
$$W=W(u) \df \inf\{v> u :
|(T(\Lambda_v) +g(v)) - (T(\Lambda_u) +g(u))|\geq 2c_2\}.$$
By the strong Markov property applied at $T(\Lambda_W)$ 
we have $P^{x_0}_h(A_W \mid W< \infty) > 1-c_2$.
Since $A_u \cap \{W< \infty\} \cap A_W = \emptyset$,
it follows that 
$P^{x_0}_h(A_W\cap\{W<\infty\})<c_2$, and hence
$P^{x_0}_h( W< \infty) < c_2/(1-c_2)$.
This proves the Lemma, since
we may assume that $c_2>0$ is arbitrarily small 
by choosing $u$ sufficiently large. \qed
\enddemo 

We now make some general observations about parabolic Martin
boundaries. Let $D$ be a domain. For 
$\phi$ a parabolic function on $\dd$, and $v<0$, define 
$$
\phi_v(x,t)\df\phi(x,t+v).
$$
Then $\phi_v$ is also parabolic. Moreover, if $\phi$ is minimal then
$\phi_v$ is either minimal or $\phi_v\equiv 0$ (see Doob (1984)
1.XV.17).

\proclaim{\Tn{159} Lemma} Let $D$ be a domain. Let $\phi$ be
parabolic on $\dd$, and let $v<0$. Then the laws of $X$ under
$P^{x,t}_{\phi_v}$ and $P^{x,t+v}_\phi$ are the same.
\endproclaim

\demo{Proof} It suffices to show that 
$P^{x,t}_{\phi_v}(A)=P^{x,t+v}_\phi(A)$, for $A$ an event of the form
$\{X(t_1)\in A_1,\dots,X(t_n)\in A_n\}$, where
$t_1<t_2<\dots <t_n$. But
$$
\align
P^{x,t}_{\phi_v}(A)&=
\frac{1}{\phi_v(x,t)}E^{x,t}[\jeden_A \phi_v(X_{t_n},\tau_{t_n})]\\
&=\frac{1}{\phi(x,t+v)}E^{x,t}[\jeden_A \phi(X_{t_n},\tau_{t_n}+v)]\\
&=\frac{1}{\phi(x,t+v)}E^{x,t+v}[\jeden_A\phi(X_{t_n},\tau_{t_n})]
=P^{x,t+v}_\phi(A). \qed
\endalign
$$
\enddemo

Now, if $(y_k,t_k)\in \dd$, $(y_k,t_k)\to z\in\prt^M\dd$, and each 
$t_k<v$, then
$$
\align
\dot K((x,t),(y_k,t_k-v))&=
\frac{p_{t-t_k+v}^D(x,y_k)}{p^D_{-t_k+v}(x_0,y_k)}\\
&=\frac{p_{t-t_k+v}^D(x,y_k)}{p^D_{-t_k}(x_0,y_k)}\cdot
\frac{p^D_{-t_k}(x_0,y_k)}{p_{-t_k+v}^D(x_0,y_k)}\\
&\to \frac{\dot K((x,t+v),z)}{\dot K((x_0,v),z)}.
\endalign
$$
Thus, provided $\dot K((x_0,v),z)>0$, it follows that $(y_k,t_k-v)$
converges in $\dd^M$ to a point $\Phi_v z\in\prt^M\dd$ with
$$
\dot K(\,\cdot\,,\Phi_v z)=
\frac{\dot K_v(\,\cdot\,,z)}{\dot K_v((x_0,0),z)}.\tag\tn{160}
$$
Of course, it may happen that $\Phi_v z=z$. Note also that 
$$
\dot K((x_0,0),\Phi_v z)=1,\tag\tn{166}
$$
so that 
$\Phi_v z$ is a minimal point, if and only if
$\dot K_v(\,\cdot\,,z)$ is a minimal function. 

It would simplify several future arguments, if the map $\Phi_v$ could
be defined for $v>0$ as well. A natural way of doing this would
be to set
$$
\phi_v(x,t)\df
\cases \phi(x,t+v),&t+v\le 0\\
\int p^D_{t+v}(x,y)\phi(y,0)dy,& t+v>0.
\endcases
$$
The obstacle to this approach is that in general, this integral need
not converge.

The following result is well known. See,
for example, Theorems C and E of  Aronson (1968).

\proclaim{\Tn{129} Lemma} Let $D$ be a domain, and let $A\subset\dd$
be compact.
\roster
\item"(i)" Let $\eps>0$ and $M<\infty$. There exists a
$\delta>0$ such that if $u$ is parabolic on $\dd$ and $u\le M$, then
$|u(z)-u(z')|<\eps$ whenever $z,z'\in A$ and $|z-z'|<\delta$.
\item"(ii)" Let $x\in D$. There exists an $M<\infty$ such that
if $u$ is parabolic on $\dd$, and $u(x,0)\le 1$, then $u\le M$ on $A$.
\endroster
\endproclaim

\proclaim{\Tn{164} Lemma} Let $D$ be a domain. 
Suppose that $(y_k,t_k)\in\dd$
converge to some $z\in\prt^M\dd$, and that $a_k\to 0$. 
Let $(x,t)\in\dd$ (so that, in particular, $t<0$) and suppose
that $\dot K((x,t),z)>0$. Then
$$
\frac{p^D_{t-t_k}(x,y_k)}{p^D_{t-t_k-a_k}(x,y_k)}\to 1\tag\tn{132}
$$
as $k\to\infty$. Moreover,
$$
(y_k,t_k+a_k-t)\to\Phi_t z. \tag\tn{128}
$$ 
\endproclaim

\demo{Proof} If $\dot K((x,t),z)>0$, then the 
$\dot K((x,t),(y_k,t_k))$ are bounded away from $0$. Since 
$\dot K((x_0,0),(y_k,t_k))=1$, (ii) of Lemma \lbl{129} shows
that the $\dot K(\,\cdot\,,(y_k,t_k))$ are uniformly bounded, on a 
suitable neighbourhood of $(x,t)$. Applying (i) of Lemma \lbl{129} on
this neighbourhood shows that
$$
\frac{p^D_{t-t_k}(x,y_k)}{p^D_{t-t_k-a_k}(x,y_k)}
=\frac{\dot K((x,t),(y_k,t_k))}{\dot K((x,t-a_k),(y_k,t_k))}\to 1,
$$
as $k\to\infty$, showing \Lbl{132}.

To prove \Lbl{128}, we must show that 
$$
\lim_{k\to\infty} \dot K((x,s),(y_k,t_k+a_k-t))
=\lim_{k\to\infty} \dot K((x,s),(y_k,t_k-t))
$$
for every $(x,s)\in\dd$. But as before,
$$
\align
\dot K((x,s),&(y_k,t_k+a_k-t))
=\frac{p_{s+t-t_k-a_k}(x,y_k)}{p_{t-t_k-a_k}(x_0,y_k)}\\
&=\frac{\dot K((x,s+t-a_k),(y_k,t_k))}{\dot K((x_0,t-a_k),(y_k,t_k))}\\
&\to \frac{\dot K((x,s+t),z)}{\dot
K((x_0,t),z)}=\dot K((x,s),\Phi_t z). \qed
\endalign
$$ 
\enddemo

\proclaim{\Tn{131} Lemma} Assume that $f(u)\to 0$ as $u\to\infty$.
Let $(z_k,t_k) \in\dd_f$ converge to 
$z\in\prt^M\dd_f$, and suppose that 
$\dot K((x,t),z)>0$ for every $(x,t)\in\dd_f$. If
$y_k\in D_f$ and $z_k^d=y_k^d$ for each $k$, then for some
$c_1<\infty$ and $c_2>0$,
and for every $q<0$ and $(x,t)\in\dd_f$,
$$
\align
\limsup_{k\to\infty} \dot K((x,t),(y_k,t_k-q))
&\le c_1 \dot K((x,t),\Phi_q z), \tag\tn{133}\\
\liminf_{k\to\infty} \dot K((x,t),(y_k,t_k-q))
&\ge c_2 \dot K((x,t),\Phi_q z).\tag\tn{134}
\endalign
$$ 
\endproclaim

\demo{Proof}
Let $r_0>0$ be so small that
for each $w\in\prt D_f$, the set $\prt D_f\cap B(w,r_0 f(w^d))$ is
the graph of a Lipschitz function $F$, with Lipschitz constant
$\lambda_0$ in some orthonormal coordinate system $\text{\sl CS}_w$.
Let the coordinates of $x$ in $\text{\sl CS}_w$ be $(\hat x,x')$, so
that
$$
D_f\cap B(w,r_0 f(w^d))=\{(\hat x,x') : x'> F(\hat x)\}
\cap B(w,r_0 f(w^d)).
$$

Let 
$$
\align
\Psi_r(w,s) &= \{ (x,t)\in \dd_f:  |x-w| < r, 
|s-t| < r^2 \}, \\
A_r(w) &= (\hat w, w' +r) \quad \text{in} \ CS_w. 
\endalign
$$
We fix a suitable $\bar s<0$ and 
apply Theorem 1.6 of Fabes et al. (1986) to some 
$\Psi_{r/8}(w,\bar s)$, to see that if 
$x_1,x_2\in D_f$, $w\in\prt D_f$, $r<r_0f(w^d)/2$,
$s,s'<\bar s$ and 
$y\in B(w,r/8)$, then
$$
\aligned
\frac{p^{D_f}_{-s}(x_2,y)}{p^{D_f}_{-s'}(x_1,y)}
&=\frac{\dot G_{\dd_f}((y,\bar s),(x_2,\bar s + s))}
{\dot G_{\dd_f}((y,\bar s),(x_1,\bar s + s'))}\\ 
&\leq c_1  \frac{\dot G_{\dd_f}((A_r(w),\bar s+2r^2),
(x_2,\bar s + s))}
{\dot G_{\dd_f}((A_r(w),\bar s-2r^2),(x_1,\bar s + s'))}\\
&=c_1 \frac{p^{D_f}_{-s+2r^2}(x_2,A_r(w))}
{p^{D_f}_{-s'-2r^2}(x_1,A_r(w))}.
\endaligned\tag\tn{153}
$$
Note that, although Theorem 1.6 of Fabes et al. (1986) would in
principle allow the above constant
$c_1$ to depend on $f(w^d)$, in fact a scaling argument shows that it
does not.

Fix $(x,t)\in\dd_f$ and $q<0$.
Let $M = \bigcup _{w\in \prt D_f} B(w, r_0 f(w^d)/32)$. If $y_k\in
M$, choose $w$ so that $y_k\in B(w, r/8)$, where $r=r_0f(w^d)/4$.
With this choice of $r$, set
$$
\bar y_k=A_r(w),\qquad a_k=2r^2.
$$
If $y_k\notin M$, set
$$
\bar y_k=y_k,\qquad a_k=0.
$$
The assumption that
$\dot K((x,t),z)>0$ for every $(x,t)\in\dd_f$
easily implies that $t_k \to - \infty$.  
By \Lbl{153},
$$
\frac{p^{D_f}_{t+q-t_k}(x,y_k)}{p^{D_f}_{q-t_k}(x_0,y_k)}
\le c_1 \frac{p^{D_f}_{t+q-t_k+a_k}(x,\bar y_k)}
{p^{D_f}_{q-t_k-a_k}(x_0,\bar y_k)}, \tag\tn{161}
$$
for $k$ so large that $t_k-t-q<\bar s$.

Let $b_k=f^2(z^d_k)$. A precise version of the parabolic Harnack
principle (see Theorem 0.2 of Fabes et al. (1986)) implies that for
$k$ large and for every $v\in D_f$ with 
$|v^d-z^d_k|<f(z^d_k)$ and $v\notin M$, we have
$$
\frac{p^{D_f}_{t+q-t_k+a_k}(x,\bar y_k)}
{p^{D_f}_{q-t_k-a_k}(x_0,\bar y_k)}
\le c_2 \frac{p^{D_f}_{t+q-t_k+a_k+b_k}(x,v)}
{p^{D_f}_{q-t_k-a_k-b_k}(x_0,v)}. \tag\tn{162}
$$

As above, take $\bar z_k$ equal to either $z_k$ (if $z_k\notin M$), or
an $A_r(w)$ (if $z_k\in B(w,r/8)$, where $r=r_0 f(w^d)/4$).
Take  $d_k$ equal to $0$ or $2r^2$ respectively.
Therefore
$$
\frac{p^{D_f}_{t+q-t_k+a_k+b_k}(x,\bar z_k)}
{p^{D_f}_{q-t_k-a_k-b_k}(x_0,\bar z_k)}
\le c_1\frac{p^{D_f}_{t+q-t_k+a_k+b_k+d_k}(x,z_k)}
{p^{D_f}_{q-t_k-a_k-b_k-d_k}(x_0,z_k)}, \tag\tn{163}
$$
for $k$ large, as before. Since $q<0$, $a_k\to 0$, $b_k\to 0$, and 
$d_k\to 0$, it follows
from \Lbl{132} that 
$$
\align
\lim_{k\to\infty}\frac{p^{D_f}_{t+q-t_k+a_k+b_k+d_k}(x,z_k)}
{p^{D_f}_{q-t_k-a_k-b_k-d_k}(x_0,z_k)}
&=\lim_{k\to\infty}\frac{p^{D_f}_{t+q-t_k}(x,z_k)}
{p^{D_f}_{q-t_k}(x_0,z_k)}\\
&=\lim_{k\to\infty}\dot K((x,t),(z_k,t_k-q))
=\dot K((x,t),\Phi_q z).
\endalign
$$
Thus, taking $v=\bar z_k$, it follows from this
and \Lbl{161}-\Lbl{163} that 
$$
\limsup_{k\to\infty} \dot K((x,t),(y_k,t_k-q))=
\limsup_{k\to\infty}
\frac{p^{D_f}_{t+q-t_k}(x,y_k)}{p^{D_f}_{q-t_k}(x_0,y_k)}
\le c_3 \dot K((x,t),\Phi_q z)
$$
as well, proving \Lbl{133}. The argument for \Lbl{134} is similar.\qed
\enddemo

We may improve upon the conclusion of Lemma \lbl{131}, 
by assuming that $z$ is minimal:

\proclaim{\Tn{135} Lemma} 
Assume that $f(u)\to 0$ as $u\to\infty$.
Let $(z_k,t_k) \in\dd_f$ converge to minimal point
$z\in\prt^M_0\dd_f$, and suppose that 
$\dot K((x,t),z)>0$ for every $(x,t)\in\dd_f$. If
$y_k\in D_f$ satisfy $z_k^d=y_k^d$ for each $k$, and $q_k\to q<0$,
then $(y_k,t_k-q_k)\to\Phi_q z$. That is,
$$
\lim_{k\to\infty} \dot K((x,t),(y_k,t_k-q_k))
= \dot K((x,t),\Phi_q z)
$$ 
for every $(x,t)\in\dd_f$.
\endproclaim

\demo{Proof} 
We first consider the limit of $(y_k,t_k-q)$. If $w$ is any limit 
point of this sequence, then by \Lbl{133} we have that
$$
\dot K(\,\cdot\,,w)\le c_1 \dot K(\,\cdot\,,\Phi_q z).
$$
By minimality of $z$ (and hence $\Phi_q z$), in fact
$$
\dot K(\,\cdot\,,w) = c \dot K(\,\cdot\,,\Phi_q z)
$$
for some $c<\infty$. By \Lbl{134} we must have $c>0$, so $w\neq\dot 0$.

Let $k_i$ be a subsequence along which
$(y_{k_i},t_{k_i}-q)\to w$. By passing to a further subsequence,
if necessary, we may also ensure that
$(y_{k_i},t_{k_i}-q/2)$ converges to some $w'\neq\dot 0$. Then
$w=\Phi_{q/2} w'$, so by \Lbl{166}, 
$$
\dot K((x_0,0),w)=1=\dot K((x_0,0),\Phi_q z).
$$
Thus $c=1$, and so $w=z$. Since $\Phi_q z$ is the only limit point
of $(y_k,t_k-q)$, it follows that the sequence itself converges to
$\Phi_q z$.

Similarly, $(y_k,t_k-q/2)\to\Phi_{q/2} z$.
Since $\Phi_q z=\Phi_{q/2}(\Phi_{q/2} z)$, we may set $a_k=q-q_k$,
and apply \Lbl{128} (with $t=q/2$), to obtain in addition that
$(y_k,t_k-q_k)\to\Phi_q z$, as required. \qed
\enddemo

\demo{Proof of Theorem \slbl{1}{14}}

(i) 
Assume either (a) or (b) of (i) of the Theorem,
and recall that this implies that $f(u)\to 0$ as $u\to\infty$.

Let $\Cal H_s$ denote the set of points $z$ of the minimal Martin
boundary $\prt^M_0\dd_f$, such that $g(u)-\tau(\Lambda_u)\to s$, 
$P^{x_0,0}_z$-a.s.
Set $\Cal H=\bigcup_{s\in\R}\Cal H_s$.
Recall that if $\phi$ is a minimal parabolic
function, then the tail
$\sigma$-field of every $\phi$-transform of space-time
Brownian motion is trivial. By Lemma \Lbl{148}, the random 
variable 
$\lim_{u\to\infty}g(u)-\tau(\Lambda_u)$ is well defined
$P^{x_0}_h$-a.s. It is clearly measurable with respect to the tail 
$\sigma$-field of $\dot X_t$, and so 
$$
h(x)=\int_{\Cal H}\dot K(\,\cdot\, ,z)\mu(dz),
$$
for some measure $\mu$ concentrated on $\Cal H$. In particular, it
follows that 
$\Cal H_s$ is non-empty, for some $s\in\R$. 
We will work towards proving that, in fact, 
$$
\text{every $\Cal H_s$ consists
of a single point,}\tag\tn{130}
$$ 
namely the $z_s$ of (A). 

In fact, the conclusion of (B) will follow immediately from \Lbl{130},
since $\Cal H_{s_1}$ and $\Cal H_{s_2}$
are disjoint if $s_1\neq s_2$. 

For $z\in\Cal H_s$, we have that 
$g(u)-\tau(\Lambda_u)\to s$, $P^{x_0,0}_z$-a.s. A standard argument
now shows that the same is true $P^{x,t}_z$-a.s., for every
$(x,t)\in\dd_f$. Thus (D) will also follow
immediately, once \Lbl{130} is proven.

(ii) 
It is a routine matter to prove that if 
$$
P^{x_0,0}_h\left(\lim_{u\to\infty} (g(u)-\tau(\Lambda_u))
\in (s_1,s_2)\right) >0\tag\tn{127}
$$
then for every $s_3 \in \R$,
$$
P^{x_0,0}_h\left(\lim_{u\to\infty} (g(u)-\tau(\Lambda_u))
\in (s_1+s_3,s_2+s_3)\right) >0.
$$
Hence, \Lbl{127} holds for all $-\infty < s_1 < s_2 < \infty$.
Therefore
$$
\mu\left(\bigcup_{s\in(s_1,s_2)}\Cal H_s\right)>0,
$$
for every such $s_1,s_2$. This will establish (C). Moreover, it
shows that 
$$
\text{$\exists\  \{s_k\}_{k\geq 1}$ such that
$\lim_{k\to \infty } s_k = \infty$ and for every $k$,
$\Cal H_{s_k}\neq\emptyset$.}\tag\tn{136}
$$

If $\phi=\dot K(\,\cdot\,,z)$, where $z\in\Cal H_s$, and $v<0$, then
by Lemma \lbl{159}, 
$$
\align
1&=P^{x,t+v}_\phi(g(u)-\tau(\Lambda_u)\to s)
=P^{x,t+v}_\phi(g(u)+T(\Lambda_u)-t-v\to s)\\
&=P^{x,t}_{\phi_v}(g(u)+T(\Lambda_u)-t-v\to s)
=P^{x,t}_{\phi_v}(g(u)-\tau(\Lambda_u)\to s+v).
\endalign
$$
That is, 
the pole of $\phi_v$ belongs to $\Cal H_{s+v}$.
Thus, $\Phi_v$ maps 
$\Cal H_s$ into
$\Cal H_{s+v}$. Appealing to \Lbl{136}, we conclude that
$\Cal H_s$ is nonempty, for every $s\in\R$.

(iii) Let $s\in\R$, and pick $z\in\Cal H_s$. 
For any sequence $u_k\to\infty$, we may set
$y_k=X(T(\Lambda_{u_k}))$, and
$t_k=\tau(\Lambda_{u_k})$.
Because $\dot X(T(\Lambda_{u_k}))\to z$
in the Martin topology, $P^{x,0}_z$-a.s., it follows that
we have constructed a sequence
$(y_k,t_k)\to z$ as in part (A), with $s_k\df g(y^d_k)-t_k\to s$.

Next we will show that ${\Cal H}_s$ consists of a single point,
for each $s$. Suppose to the contrary that $z, \wt z \in{\Cal H}_s$
for some $s$. It is easy to see that we must have
$\Phi_{s'-s} z \ne \Phi_{s'-s} \wt z$ for some $s'<s$.
Fix any sequence $u_k\to\infty$.
Consider any sequence $(y_k,t_k)\to z$, with $ g(y^d_k)-t_k\to s$
and $y_k^d = u_k$,
constructed as in the previous paragraph. Let
$(\wt y_k,\wt t_k)$ be the analogous sequence with
$(\wt y_k,\wt t_k)\to \wt z$, $\wt y_k^d = u_k$ and $ g(\wt y^d_k)
-\wt t_k\to s$.
Note that $t_k - \wt t_k \to 0$ because
$y_k^d = \wt y_k^d = u_k$, $ g(y^d_k)-t_k\to s$,
and $ g(\wt y^d_k)-\wt t_k\to s$. Lemma \lbl{135}
implies that $(y_k, t_k + s - s') \to \Phi_{s'-s} z$.
But it also implies that 
$(y_k, t_k + s - s') \to \Phi_{s'-s} \wt z$, because
$y_k^d = \wt y_k^d$ and $t_k - \wt t_k \to 0$.
This contradicts the fact that $\Phi_{s'-s} z \ne \Phi_{s'-s} \wt z$
and so it proves our claim, and establishes \Lbl{130}.

Now let $r_k\to s$, and consider any sequence $x_k$ such that
$x_k^d\to\infty$. 
Our goal is to show that $(x_k,g(x_k^d)-r_k)\to z$, where $z$
is the only element of ${\Cal H}_s$.
Set $u_k=x_k^d$, and this time choose $s'>s$. Let
$z'$ be the element of $\Cal H_{s'}$. Note that $\Phi_{s-s'} z' = z$.
By the argument of the first paragraph
of (iii),
we may choose $(y_k,t_k')\to z'$ with $y_k^d=u_k=x_k^d$
and 
$ g(y^d_k)-t_k'\to s'$. 
Since $t_k' - [g(x_k^d)-r_k] \to -s' + s$, 
we may apply Lemma \lbl{135} and obtain that
$$
(x_k,g(x_k^d)-r_k)\to \Phi_{s-s'} z' = z.
$$ 
This finishes the proof of (A). Thus, part (i) of Theorem
\slbl{1}{14} is proven.

(iv) Turning to part (ii) of Theorem \slbl{1}{14}, suppose that
$\int_u^\infty f^3(v)dv = \infty$ for all
$u<\infty$. We also assume, as it simplifies the proof, that
$f(u)\to 0$ as $u\to\infty$. At the end we will sketch out how to
extend the argument to the general case, that
$\limsup_{u\to\infty}f(u)<\infty$.

We use a coupling argument. Fix $x_1,x_2\in D_f$, and $s\le 0$.
Let $X_1$ and $X_2$ be
independent processes, under a probability measure $P$, with the same
distributions as
$X$ under $P^{x_1,s}_h$ and $P^{x_2,s}_h$ respectively.
Thus, $\dot X_1(t)=(X_1(t),\tau_t)$ and $\dot
X_2(t)=(X_2(t),\tau_t)$ are versions of $\dot X$, where
$\tau(t)=s-t$.
Define
$$
W=\inf\{t>0 : 
X^d_1(t)=X^d_2(t)\}.
$$ 
We will show that
$$
P(W<\infty)=1.\tag\tn{146}
$$

Write $T_j(\Lambda_u)$ for the hitting time of $\Lambda_u$ by $X_j$.
We may assume, without loss of generality, that 
$x_1^d \leq x_2^d$. Set 
$u_0 = x_2^d + f(x_2^d)$,
$Y_j = T_j(\Lambda_{u_0})$ and 
$Z_j = T_j(\Lambda_u)-T_j(\Lambda_{u_0})$, where the value of
$u$ will be chosen later. 
A standard application of the boundary
Harnack principle \slbl{2}{81} shows that the Radon-Nikodym derivative
of the hitting distributions of $\Lambda_{u_0}$ under $P^{y_1}_h$
and $P^{y_2}_h$ is bounded below by $c_1>0$ for all
$y_1,y_2 \in \Lambda_{x_2^d}$. 

Let $c_2$ be so large that
$$P(Y_1 - Y_2\geq c_2 ) < c_1 /16.\tag\tn{155}$$
Use Theorem \SLbl{1}{12} (v) to find $u$ so large that for every
$v\in\R$  we have
$$P(Z_2 \in (v,v+ c_2)) < c_1/8.\tag\tn{156}$$
Let $v_1$ be the median of $Z_1$, in other words,
$$P(Z_1 \leq v_1) \geq 1/2, \quad P(Z_1 \geq v_1) \geq
1/2.\tag\tn{157}$$ 
By applying the strong Markov property at
$T(\Lambda_{u_0})$, and by our choice of $c_1$, we have 
$P(Z_2 \geq v_1) \geq c_1/2$. Now we use
\Lbl{156} to obtain that
$P(Z_2 \geq v_1+c_2) \geq 3c_1/8$. This, \Lbl{157} and
the independence of $Z_1$ and $Z_2$ show that
$$
P(Z_2-Z_1\ge c_2)\ge
P(Z_1 \leq v_1, Z_2 \geq v_1+c_2) \geq 3c_1/16.
$$
Inequality \Lbl{155} now implies that
$$
\aligned
P(T_1(\Lambda_u)<T_2(\Lambda_u))&=
P(Y_1 + Z_1 <Y_2 + Z_2)\\
&\ge P(Y_1-Y_2<c_2\le Z_2-Z_1)\\
&\ge P(Z_2-Z_1\ge c_2) -  P(Y_1-Y_2\ge c_2)\ge c_1/8.
\endaligned\tag\tn{158}
$$

Let $V_j^0 = x_j$, $\tau^0 = s$,
$T^1 = \max(T_1(\Lambda_u),T_2(\Lambda_u))$, $\tau^1=\tau(T^1)$,
$V_j^1 = X(T_j(\Lambda_u))$, $U^1=u$.
Repeat the above argument, starting from $(V_j^1,\tau^1)$ in
place of $(V_j^0,\tau^0)$, and ensuring that $U^2$ is chosen
so large that each $T_j(\Lambda_{U^2})>T^1$.
Then continue this procedure
inductively, to obtain sequences of random variables $V_j^k$,
$T^k$, $\tau^k$, and $U^k$. By the strong Markov property, \Lbl{158}
becomes that
$$
P(T_1(\Lambda_{U^{k+1}})<T_2(\Lambda_{U^{k+1}})
\mid
\Cal F_{T^k})\ge c_1/8,
$$
where $\Cal F_t$ is the filtration of $(X_1(t),X_2(t))$.
It follows that an infinite number of these events will occur,
$P$-a.s. The same is true when the roles of $X_1$ and $X_2$ are
reversed. Thus \Lbl{146} holds.

(v) According to \Lbl{146}, used repeatedly, there are
points $(x_{j,k},t_k)$ on the paths of $\dot X_j$ such that 
$x^d_{1,k}=x^d_{2,k}\to\infty$. 
Using Lemma \lbl{135}, 
as in the argument of section (iii) above, we get that
$(x_{1,k},t_k)$ and $(x_{2,k},t_k)$ have the same limit in 
$\prt^M_0\dd_f$. Thus, the limits of $\dot X_1(t)$ and $\dot X_2(t)$ 
in $\prt^M_0\dd_f$, as $t\to\infty$, are the same. Since 
$\dot X_1$ and $\dot X_2$ are independent, the measure $\mu$
such that
$h(x)=\int_{\prt^M_0 \dd_f}\dot K((x,0),z)\mu(dz)$ must actually
be supported on a singleton. That is, $h$ must be minimal as a
parabolic function.

It is the use of Lemma \lbl{135} 
that requires the assumption that
$f(u)\to 0$. If only $\limsup_{u\to\infty}f(u)<\infty$, we modify
the argument as follows. For any $\eps>0$,
$$
P_h^{x,t}(X(f((x^d)^2))\in
B((0,\dots,0,x+f(x^d)),\eps f(x^d)))\ge c(\eps)>0,\tag\tn{147}
$$
for every $(x,t)\in\dd_f$. 
Let
$$
W_\eps\df \inf\{t>0 : 
\text{$X^d_1(t),X^d_2(t)\in B((0,\dots,0,u+f(u)),\eps f(u))$
for some $u$}\}.
$$
Applying \Lbl{147} to $x=X_j(W)$ and using
another iterative argument, one can show that
$P(W_\eps<\infty)=1$
for every $\eps>0$. Taking a sequence $\eps_k\to 0$, this now
gives sequences 
$(x_{j,k},t_k)$ on the paths of $\dot X_j$, such that 
$$
x_{j,k}\in B((0,\dots,0,u_k+f(u_k), \eps_k f(u_k)),
$$
where
$u_k\to\infty$. An argument as in the proof of 
Lemmas \lbl{131} and \lbl{135} 
now shows that the
$(x_{1,k},t_k)$ and $(x_{2,k},t_k)$ have the same limit in 
$\prt^M_0\dd_f$. As before, this shows that $h$ is parabolically
minimal. \qed

\enddemo

\enddocument